\documentclass[11pt,a4paper]{article}
\pdfoutput=1
\RequirePackage{amsmath}%
\usepackage{amsfonts}%
\usepackage{amssymb}%
\usepackage{graphicx}%
\usepackage{url}%
\usepackage{color}

\newtheorem{example}{Example}

\newtheorem{remark}{Remark}

\usepackage[autolinebreaks]{mcode}
\usepackage{url}
\usepackage[title]{appendix}%

\DeclareMathOperator{\capa}{\mathrm{cap}}
\let\Re\relax
\DeclareMathOperator{\Re}{\mathrm{Re}}
\let\Im\relax
\DeclareMathOperator{\Im}{\mathrm{Im}}
\DeclareMathOperator{\diam}{\mathrm{diam}}

\renewcommand{\i}{\mathrm{i}}
\newcommand{\bI}{{\bf I}}
\newcommand{\bM}{{\bf M}}
\newcommand{\bN}{{\bf N}}
\newcommand{\CC}{{\mathbb C}}
\newcommand{\RR}{{\mathbb R}}
\newcommand{\DD}{{\mathbb D}}

\newcommand{\bt}{{\bf t}}
\newcommand{\bs}{{\bf s}}
\newcommand{\bx}{{\bf x}}
\newcommand{\br}{{\bf r}}

\begin{document}
\title{Fast computation of analytic capacity}
	\author{Mohamed M S Nasser$^{\rm a}$, Christopher C. Green$^{\rm a}$, Matti Vuorinen$^{\rm b}$}
	
	\date{}
	\maketitle
 	
	\vskip-0.8cm %
	\centerline{$^{\rm a}$Department of Mathematics, Statistics \& Physics, Wichita State University,} %
	\centerline{Wichita, KS 67260-0033, USA}%
	\centerline{\tt mms.nasser@wichita.edu, christopher.green@wichita.edu}%
	\centerline{$^{\rm b}$Department of Mathematics and Statistics, University of Turku, Turku, Finland} %
	\centerline{\tt vuorinen@utu.fi}%

\begin{abstract}
A boundary integral equation method is presented for fast computation of the analytic capacities of compact sets in the complex plane. The method is based on using the Kerzman--Stein integral equation to compute the Szeg\"o kernel and then the value of the derivative of the Ahlfors map at the point at infinity. The proposed method can be used for domains with smooth and piecewise smooth boundaries. When combined with conformal mappings, the method can be used for compact slit sets. Several numerical examples are presented to demonstrate the efficiency of the proposed method. We recover some known exact results and corroborate the conjectural subadditivity property of analytic capacity.
\end{abstract}

\begin{center}
\begin{quotation}
\end{quotation}
\end{center}

\section{Introduction} \label{sec:int}

Capacities -- such as analytic, logarithmic, and conformal -- are fundamental tools in complex analysis and have several applications to problems in different fields, e.g., in approximation theory, potential theory, electronics, and fluid dynamics~\cite{Crowdy1,Dav,Gar,Mur94,Tol,Za}. 
These capacities can be expressed explicitly in only a handful of special cases, and therefore numerical methods are needed to compute these capacities in the majority of instances. 

The Riemann mapping theorem states that any unbounded simply connected domain $G$ in the extended complex plane $\widehat{\CC}=\CC\cup\{\infty\}$ with $\infty\in G$ and whose boundary consists of more than one point can be mapped one-to-one onto the unit disk $\DD$ by a conformal map $f$. If we assume that
\begin{equation}\label{eq:cond}
	f(\infty)=0 \quad {\rm and}\quad f'(\infty)>0,
\end{equation}
then this mapping $f$ is unique and known as the Riemann mapping function. Here, the derivative of analytic function $f$ at the point at infinity is 
\begin{equation}\label{eq:f'inf}
	f'(\infty)=\lim_{z\to\infty}z[f(z)-f(\infty)].
\end{equation}
The so-called Ahlfors map can be regarded as an extension of the Riemann mapping function for multiply connected domains. That is, given an unbounded multiply connected domain $G$ of connectivity $m$, the associated Ahlfors map is the unique analytic function function $f$ that maps $G$ onto $\DD$ such that~\cite{Gar,Tol,Za}
\begin{equation}\label{eq:cond-A}
	f(\infty)=0,\quad f'(\infty)>0, \quad \text{and\; $f'(\infty)$ is maximal}.
\end{equation}
The Ahlfors map $w=f(z)$ is then an $m$-to-one covering of $G$ onto $\DD$ which maps each boundary component of $G$ one-to-one onto the unit circle.

Let us now introduce the notion of analytic capacity. Let $E$ be a compact subset of the complex plane $\CC$  and let $G = \widehat{\CC} \backslash E$ be  its complement in the extended complex plane $\widehat{\CC}=\CC\cup\{\infty\}$. The analytic capacity of $E$ is defined to be~\cite[p.~15]{Tol}
\[
\gamma(E)=\sup|f'(\infty)|
\]
where the supremum is taken over all analytic functions $f:G\to \CC$ such that $|f(z)|\le1$ for all $z\in G$, and $f'(\infty)$ is as in~\eqref{eq:f'inf}. In this paper, we assume that $G$ is an unbounded multiply connected domain of connectivity $m$ with $\infty\in G$.

It is well-known that analytic capacity is inextricably linked to the Ahlfors map~\cite{Dav,Gar,Mur94,Tol,Za}. If $w=f(z)$ is the unique Ahlfors map from the unbounded multiply connected domain $G$ in the $z$-plane onto the unit disk $\DD$ in the $w$-plane satisfying the conditions~\eqref{eq:cond-A}, then the analytic capacity of the set $E$ is given by~\cite{Gar,Tol} 
\begin{equation}\label{eq:cap-def}
	\gamma(E) = f'(\infty).
\end{equation}
In particular, when $E$ is compact and connected such that $G = \widehat{\CC} \backslash E$ is a simply connected domain, the Ahlfors map $w=f(z)$ is a conformal map from $G$ onto the unit disk $\DD$ and hence the analytic capacity $\gamma(E)$ is equal to the logarithmic capacity $\capa(E)$ of $E$.
However, for a general compact set $E$, we have $\gamma(E)\le \capa(E)$. See~\cite[p.~13]{Za} for details.

Closed-form expressions for analytic capacity are special and known only in a handful of cases, and it is informative to survey some of these briefly here. If $E$ is a disk of radius $r$, then~\cite[p.~17]{Tol}
\begin{equation}\label{eq:cap-cir}
	\gamma(E)=r, 
\end{equation}
and if $E$ is a square with sides of length $\ell$, then~\cite{LSN17}
\begin{equation}\label{eq:cap-sq}
	\gamma(E)=\frac{\ell\,\Gamma^2(1/4)}{4\sqrt{\pi^3}},
\end{equation}
where $\Gamma(\cdot)$ is the gamma function. 
For a complex line segment $E=[a,b]\subset\CC$, we have~\cite[p.~17]{Tol}
\begin{equation}\label{eq:cap-seg}
	\gamma(E)=\frac{1}{4}|b-a|.
\end{equation}
For any compact subset $E$ of $\RR$, we have \cite[Chapter I, Theorem 6.2]{Gar}
\[
\gamma(E)=\frac{1}{4}|E|.
\]
In particular, for $m$ non-overlapping real intervals $E_j=[a_j,b_j]$ with $a_1<b_1<\cdots<a_m<b_m$, if $E=\bigcup_{j=1}^mE_j$, we have~\cite[p.~14]{Za}: 
\begin{equation}\label{eq:cap-int}
	\gamma(E)=\frac{1}{4}|E|=\frac{1}{4}\sum_{j=1}^m|E_j|=\frac{1}{4}\sum_{j=1}^m(b_j-a_j).
\end{equation}	
If $E\subset F\subset\CC$, then $\gamma(E)\le \gamma(F)$. Furthermore, for all $z,\lambda\in\CC$, $\gamma(z+\lambda E)= |\lambda| \gamma(E)$.
For a compact and connected set $E\subset\CC$, we have
\[
\diam(E)/4\le\gamma(E)\le \diam(E)
\]
where $\diam(E)$ denotes the diameter of $E$. For more details, see~\cite[p.~17]{Tol}.

When $E$ and $F$ are disjoint connected compact subsets of $\CC$, Suita~\cite{Suita} proved that the subadditivity property
\begin{equation}\label{eq:sub}
	\gamma(E\cup F)\le \gamma(E)+\gamma(F)
\end{equation}
holds. For general compact sets, the proof of this property is still an open problem. 
However, Tolsa~\cite{Tol03} proved the semi-additivity of analytic capacity: there exists a constant $c$ such that for all compact
sets $E,F\subset\CC$, the analytic capacity satisfies
\begin{equation}\label{eq:subc}
	\gamma(E\cup F)\le c\left(\gamma(E)+\gamma(F)\right).
\end{equation}
Moreover, proving the conjectural subadditivity of analytic capacity for arbitrary compact
sets $E,F\subset\CC$ is equivalent to proving it for all disjoint compact sets $E,F\subset\CC$ that are finite unions of disjoint closed disks, all with the same radius~\cite{Mel,YR13}. Several numerical examples have been considered by Younsi \& Ransford in~\cite{YR13} for purely circular compact sets. Numerical results for sets other than circular ones were also presented in~\cite{YR13,YR18}.
All of these examples provide convincing evidence to suggest that the conjectural subadditivity property for analytic capacity is true. It should also be pointed out that, from~\eqref{eq:cap-int}, in the case of multiple real slits, equality as opposed to inequality holds in \eqref{eq:sub}.
The subadditivity property will be of significant interest in the ensuing discussion and will be corroborated numerically in several cases.

Several numerical methods are available in the literature for computing the logarithmic and conformal capacities. One of these methods is based on the boundary integral equation (BIE) with the generalized Neumann kernel~\cite{LSN17,Nvm}. For the numerical computation of analytic capacity, to the best of our knowledge, the only available numerical method is that given in~\cite{YR13}. The method is based on using quadratic minimization for the numerical computation of upper and lower bounds for the analytic capacity, which, in principle, converges to its exact value.
This method has been used before to compute the logarithmic capacity~\cite{Ran10}.
In this paper, we present a fast and accurate BIE method for the numerical computation of the analytic capacity. Our presented method will be used to compute the analytic capacity for a wider class of compact sets, including those with smooth boundaries, piecewise smooth boundaries, and sets consisting of only slits. The method is based on using a BIE for the Szeg\"o kernel (refer to~\cite{BelBook,Gara,Ker-St,Ker-Tru,Mur94} for the definition and basic properties of the Szeg\"o kernel). The BIE has been used by Kerzman \& Trummer~\cite{Ker-Tru} to compute the conformal mapping for simply connected domains. Bell~\cite{BelAhl} proved that the BIE can be used to compute the Ahlfors map of bounded multiply connected domains. In~\cite{Visual}, the BIE has been implemented numerically to compute both the Riemann mapping function and the Ahlfors map for bounded simply and multiply connected domains, respectively. See also~\cite{BelBook,Ker-St,OD,Tru-Sze}.

Besides this introductory section, our paper is structured in the following way. 
In Section 2, we introduce a numerical method for computing the analytic capacity for compact sets bounded by smooth or piecewise smooth Jordan curves, and several numerical examples for such sets are presented in Section 3. In Section 4, we consider compact slit sets. We provide our concluding remarks in Section 5. 
Finally, in Appendix A, for a given multiply connected rectilinear slit domain, we review an iterative method from~\cite{NG18} for the construction of a preimage unbounded multiply connected domain bounded by smooth Jordan curves.

\section{The numerical method} \label{sec:nm}

\subsection{Analytic capacity and the Szeg\"o kernel} \label{sec:cap}

Let $E$ be a compact subset of the complex plane $\CC$ and let $G = \widehat{\CC} \backslash E$ be its complement  in the extended complex plane $\widehat{\CC}$. 
We assume that $G$ is an unbounded multiply connected domain bounded by $m$ smooth, or piecewise smooth, Jordan curves $\Gamma_1,\ldots,\Gamma_m$. Domains bounded by slits will be considered in Section~\ref{sec:Slits} below.

From~\eqref{eq:cap-def}, the analytic capacity $\gamma(E)$ is calculated by computing the derivative $f'(\infty)$ of the Ahlfors map $f$ from the domain $G$ onto the unit disk $\DD$ with the normalization~\eqref{eq:cond}. A BIE method for computing the Ahlfors map for a \emph{bounded} multiply connected domain is  presented in~\cite{BelAhl,Visual}; however, the domain $G$ is \emph{unbounded}, and therefore a preliminary step is required.
We first conformally map the unbounded multiply connected domain $G$ onto a bounded multiply connected domain $D$ using the M\"obius transformation
\[
\zeta=M(z)=\frac{1}{z-\alpha},
\]
where $\alpha$ is a point in the interior of any of the curves $\Gamma_j$. 
The point at infinity is mapped onto the origin, $M(\infty)=0\in D$.
Let $w=F(\zeta)$ be the Ahlfors map from the bounded multiply connected domain $D$ onto the unit disk $\DD$ such that $F(0)=0$, $F'(0)>0$, and $F'(0)$ is maximal. It follows immediately that
\[
f(z)=F(M(z)), \quad z\in G,
\]
is the Ahlfors map from the unbounded domain $G$ onto the unit disk $\DD$ which satisfies the conditions~\eqref{eq:cond-A}. 
Note that 
\[
f(\infty)=F(M(\infty))=F(0)=0
\]
and
\[
f'(\infty)=\lim_{z\to\infty}z[f(z)-f(\infty)]=\lim_{z\to\infty}zf(z)=\lim_{z\to\infty}zF(M(z)).
\]
Note also that
$\zeta=M(z)=1/(z-\alpha)$ if and only if $z=M^{-1}(\zeta)=a+1/\zeta$ and hence $z\to\infty$ if and only if $\zeta\to 0$. Thus, since $F(0)=0$,
\[
f'(\infty)
=\lim_{\zeta\to 0}\left(a+\frac{1}{\zeta}\right)F(\zeta)
=F'(0)>0.
\]
It is straightforward to prove that $f'(\infty)$ is maximal since $F'(0)$ is maximal.

As $D$ is a bounded multiply connected domain, the Ahlfors map  $w=F(\zeta)$ from $D$ onto the unit disk $\DD$ can be computed using the method presented by Bell~\cite{BelAhl} (see also~\cite{Visual}). However, computing the analytic capacity $\gamma(E)$ requires only computing the derivative $F'(0)$ since $\gamma(E)=f'(\infty)=F'(0)$. In fact, it follows from~\cite{BelAhl} that
\[
F'(0)=2\pi S(0,0)
\]
where $S(\zeta,0)$ is the Szeg\"o kernel for the bounded multiply connected domain $D$ with respect to the base point $0\in D$.
Our numerical method is based on computing $S(0,0)$ using the BIE related to the Szeg\"o kernel in multiply connected domains~\cite{BelAhl}. Then
\begin{equation}\label{eq:cap-b}
	\gamma(E) = f'(\infty)=F'(0)= 2\pi S(0,0).
\end{equation}

\subsection{Integral equation for the Szeg\"o kernel} 

Assume that each boundary component $\Gamma_j$ is parametrized by a $2\pi$-periodic function $\zeta_j(t)$, $t\in J_j=[0,2\pi]$, $j=1,\ldots,m$. For domains with corners, the parametrization $\zeta_j(t)$ is defined as described in~\cite{LSN17}.
We define the total parameter domain $J$ as the disjoint union of the $m$ intervals $J_j=[0,2\pi]$, $j=1,\ldots,m$. The whole boundary $\Gamma$ is therefore parametrized by
\begin{equation}\label{eq:zeta}
	\zeta(t)= \left\{ \begin{array}{l@{\hspace{0.5cm}}l}
		\zeta_1(t),&t\in J_1, \\
		\quad\vdots & \\
		\zeta_m(t),&t\in J_m.
	\end{array}
	\right.
\end{equation}
See~\cite{LSN17,Nas-ETNA,NG18} for more details. Further, the boundary $\partial D$ of the bounded multiply connected domain $D$ is parametrized by 
\begin{equation}\label{eq:eta}
	\eta(t)=\frac{1}{\zeta(t)-\alpha}, \quad t\in J.
\end{equation}

The Szeg\"o kernel for the bounded multiply connected domain $D$ with respect to the base point $0\in D$ can be computed by solving the second-kind Fredholm integral equation~\cite{BelAhl,Visual}
\begin{equation}\label{eq:ie-1n}
	S(\eta(t),0)+\int_{J} A(\eta(t),\eta(s))S(\eta(s),0)|\eta'(s)|ds=
	\overline{\frac{1}{2\pi\i}\frac{\eta'(t)}{|\eta'(t)|\eta(t)}}
\end{equation}
where
\[
A(\eta(t),\eta(s))=\overline{\frac{1}{2\pi\i}\frac{\eta'(t)}{|\eta'(t)|(\eta(t)-\eta(s))}}
-\frac{1}{2\pi\i}\frac{\eta'(s)}{|\eta'(s)|(\eta(s)-\eta(t))}
\]
The kernel $A(\eta(t),\eta(s))$ is continuous with $A(\eta(t),\eta(t))=0$. The integral equation in~\eqref{eq:ie-1n} is known as the Kerzman--Stein BIE~\cite{Ker-St,Ker-Tru,OD,Tru-Sze}.
It was proved in~\cite{BelAhl} that this BIE can be used also for the computation of the Szeg\"o kernel for bounded multiply connected domains.
Multiplying both sides of~\eqref{eq:ie-1n} by $\eta'(t)/\eta(t)$ and defining
\begin{equation}\label{eq:S-phi}
	\phi(t)=\frac{S(\eta(t),0)\eta'(t)}{\eta(t)},
\end{equation}
the BIE~\eqref{eq:ie-1n} can be written as
\begin{equation}\label{eq:ie-3n}
	\phi(t)
	+\int_{J} A(\eta(t),\eta(s))\frac{\eta'(t)}{\eta(t)}\frac{\eta(s)}{\eta'(s)}
	\phi(s)|\eta'(s)|ds=
	\frac{\i}{2\pi}\frac{|\eta'(t)|}{|\eta(t)|^2}
\end{equation}
Now, since
\[
\eta(t)=\frac{1}{\zeta(t)-\alpha}, \quad
\zeta(t)=\frac{1}{\eta(t)}+\alpha, 
\]
it follows that
\[
\phi(t)
+\int_{J} A\left(\frac{1}{\zeta(t)-\alpha},\frac{1}{\zeta(s)-\alpha}\right)\frac{\zeta(s)-\alpha}{\zeta(t)-\alpha}\frac{\zeta'(t)}{\zeta'(s)}
\frac{|\zeta'(s)|}{|\zeta(s)-\alpha|^2}\phi(s)ds=
\frac{\i}{2\pi}|\zeta'(t)|,
\]
or equivalently, 
\[
\phi(t)-\int_{J} \left(\overline{\frac{1}{2\pi\i}\frac{\zeta'(s)}{|\zeta'(s)|(\zeta(s)-\zeta(t))}}
-\frac{1}{2\pi\i}\frac{\zeta'(t)}{|\zeta'(t)|(\zeta(t)-\zeta(s))}\right)|\zeta'(t)|\phi(s)ds=\frac{\i}{2\pi}|\zeta'(t)|.
\]	
Taking the conjugate of both sides and then multiplying by $\i/|\zeta'(t)|$, we obtain
\[
\frac{\i\,\overline{\phi(t)}}{|\zeta'(t)|}+\int_{J} \left(\overline{\frac{1}{2\pi\i}\frac{\zeta'(t)}{|\zeta'(t)|(\zeta(t)-\zeta(s))}}
-\frac{1}{2\pi\i}\frac{\zeta'(s)}{|\zeta'(s)|(\zeta(s)-\zeta(t))}\right)\frac{\i\,\overline{\phi(s)}}{|\zeta'(s)|}|\zeta'(s)|ds=\frac{1}{2\pi}.
\]	
which can be written in the concise form
\begin{equation}\label{eq:ie-psi}
	\psi(t)+\int_{J} A(\zeta(t),\zeta(s))\psi(s)|\zeta'(s)|ds=\frac{1}{2\pi}, \qquad
	\psi(t)=\frac{\i\,\overline{\phi(t)}}{|\zeta'(t)|}.
\end{equation}
This BIE~\eqref{eq:ie-psi} is a modification of the Kerzman--Stein BIE~\eqref{eq:ie-1n}. 

By~\eqref{eq:cap-b}, computing the analytic capacity $\gamma(E)$ requires computing the value of the Szeg\"o kernel $S(0,0)$. Since the Szeg\"o kernel $S(\zeta,0)$ is an analytic function in the domain $D$, by the Cauchy integral formula, we have
\[
S(0,0)=\frac{1}{2\pi\i}\int_{\partial D} \frac{S(\zeta,0)}{\zeta}d\zeta
=\frac{1}{2\pi\i}\int_{J} \frac{S(\eta(t),0)}{\eta(t)}\eta'(t)dt.
\] 
Then, by solving the BIE~\eqref{eq:ie-psi} for $\psi(t)$ and using~\eqref{eq:S-phi} and~\eqref{eq:ie-psi}, we have
\[
S(0,0)
= \frac{1}{2\pi\i}\int_{J} \phi(t)dt
= \frac{1}{2\pi}\int_{J} \overline{\psi(t)}|\zeta'(t)|dt.
\]
It then follows at once from~\eqref{eq:cap-b} that 
\begin{equation}\label{eq:cap-mun}
	\gamma(E) = 2\pi S(0,0) = \int_{J} \overline{\psi(t)}|\zeta'(t)|dt.
\end{equation}
It is immediate from~\eqref{eq:cap-mun} that the integral $\int_{J} \overline{\psi(t)}|\zeta'(t)|dt$ must be real and hence
\begin{equation}\label{eq:cap-mu2n}
	\gamma(E) = \Re\left[\int_{J} \overline{\psi(t)}|\zeta'(t)|dt\right]
	= \int_{J} \Re[\psi(t)]|\zeta'(t)|dt.
\end{equation}

\subsection{Numerical solution of the integral equation} 

The Kerzman–Stein BIE~\eqref{eq:ie-1n} has been used to compute conformal mappings for bounded and unbounded simply connected domains~\cite{Ker-Tru,Murid1,OD,Tru-Sze}, and in~\cite{BelAhl,Visual} to compute the Ahlfors map for bounded multiply connected domains.
A combination of the usage of the Kerzman--Stein BIE~\eqref{eq:ie-1n} and the Fast Multipole Method (FMM)~\cite{Gre-Gim12,GR} has been presented in~\cite{OD} for computing conformal mappings for bounded simply connected domains.

In this paper, to compute the analytic capacity $\gamma(E)$, we will solve the BIE~\eqref{eq:ie-psi} which is a modified version of the Kerzman--Stein BIE~\eqref{eq:ie-1n}. To accomplish this, we too shall employ the FFM when solving~\eqref{eq:ie-psi}. Since the integrand in~\eqref{eq:ie-psi} is $2\pi$-periodic, the BIE~\eqref{eq:ie-psi} can be best discretized by the Nystr\"om method with the trapezoidal rule~\cite{Ker-Tru,Murid1,OD,Tru-Sze}. 

For domains with smooth boundaries, we use the trapezoidal rule with equidistant nodes. 
We discretize each interval $J_p=[0,2\pi]$, for $p=1,2,\ldots,m$, by $n$ equidistant nodes $s_1, \ldots, s_n$ where
\begin{equation}\label{eq:s_i}
	s_q = (q-1) \frac{2 \pi}{n}, \quad q = 1, \ldots, n,
\end{equation}
and $n$ is an even integer. Writing $\bs=[s_1,\ldots, s_n]$, we  discretize the parameter domain $J$ by the vector $\bt = [\bs, \bs, \ldots, \bs]$ which consists of $m$ copies of $\bs$, i.e.,
\[
\bt=[t_1,t_2,\ldots,t_{mn}]
\]
where for $p=1,2,\ldots,m$ and $q = 1, \ldots, n$,
\[
t_{(p-1)n+q}=s_q.
\]
For a real or a complex function $\mu(\zeta(t))$ defined on the boundary $\Gamma$, the trapezoidal rule then yields
\begin{equation}\label{eq:trap}
	\int_J\mu(\zeta(t))dt = \sum_{p=1}^{m}\int_{J_p}\mu(\zeta_p(t))dt
	\approx \sum_{p=1}^{m}\sum_{q=1}^{n}\frac{2\pi}{n}\mu(\zeta_p(s_q))
	=\sum_{j=1}^{mn}\frac{2\pi}{n}\mu(\zeta(t_j)).
\end{equation}
Discretizing the BIE~\eqref{eq:ie-psi} using the trapezoidal rule~\eqref{eq:trap} and substituting $t=t_i$, we obtain the linear system
\[
\psi_n(t_i)+\frac{2\pi}{n}\sum_{j=1}^{mn} A(\zeta(t_i),\zeta(t_j))|\zeta'(t_j)|\psi_n(t_j)=\frac{1}{2\pi}, \quad i=1,2,\ldots,mn,
\]
where $\psi_n$ is an approximation of $\psi$. Recall that $A(\zeta(t_i),\zeta(t_j))=0$ when $i=j$. Using the definition of the kernel $A(\zeta(t),\zeta(s))$, we have for $i=1,2,\ldots,n$,
\[
\psi_n(t_i)+\frac{\i}{n}\sum_{\begin{subarray}{c} j=1\\j\ne i\end{subarray}}^{mn} \left(\frac{\overline{\zeta'(t_i)}}{|\zeta'(t_i)|(\overline{\zeta(t_i)}-\overline{\zeta(t_j)})}
+\frac{\zeta'(t_j)}{|\zeta'(t_j)|(\zeta(t_j)-\zeta(t_i))}\right)|\zeta'(t_j)|\psi_n(t_j)=\frac{1}{2\pi}, 
\]
or equivalently,
\[
\psi_n(t_i)+\frac{\i}{n} \frac{\overline{\zeta'(t_i)}}{|\zeta'(t_i)|}\overline{\sum_{\begin{subarray}{c} j=1\\j\ne i\end{subarray}}^{mn} \frac{1}{\zeta(t_i)-\zeta(t_j)}
	|\zeta'(t_j)|\overline{\psi_n(t_j)}}
- \frac{\i}{n}\sum_{\begin{subarray}{c} j=1\\j\ne i\end{subarray}}^{mn} \frac{1}{\zeta(t_i)-\zeta(t_j)}\zeta'(t_j)\psi_n(t_j)=\frac{1}{2\pi}, 
\]
which can be written in the following concise form:
\begin{equation}\label{eq:sys-E}
	\bx+\frac{\i}{n} (\overline{\zeta'(\bt)}./|\zeta'(\bt)|).*\overline{B\left(
		|\zeta'(\bt)|.*\overline{\bx}\right)}
	- \frac{\i}{n}B\left(
	\zeta'(\bt).*\bx\right)=\br.
\end{equation}
Here, $.*$ and $./$ are the MATLAB element-by-element multiplication and division, respectively, $\bx=\psi_n(\bt)$, $\br$ is the $mn\times 1$ vector with entries $\br_{i}= \frac{1}{2\pi}$, and $B$ is the $mn\times mn$ matrix with entries 
\begin{equation}\label{e:fmm-E}
	(B)_{ij}= \left\{
	\begin{array}{l@{\hspace{1cm}}l}
		\displaystyle 0,                             & i=j,     \\[0.00cm]
		\displaystyle\frac{1}{\zeta(t_i)-\zeta(t_j)}, & i\ne j, \quad i,j=1,2,\ldots,(m+1)n.   \\
	\end{array}\right.
\end{equation}

The linear system~\eqref{eq:sys-E} will be solved using the GMRES iterative method~\cite{gmres} where the matrix-vector product can be computed using the FMM. The FMM allows us to use the method presented to compute the analytic capacity of compact sets consisting of a high number of components.
If we define the left-hand side of~\eqref{eq:sys-E} to be a function of the unknown vector $\bx$,  
\[
\mathcal{F}(\bx )=\bx+\frac{\i}{n} (\overline{\zeta'(\bt)}./|\zeta'(\bt)|).*\overline{B\left(
	|\zeta'(\bt)|.*\overline{\bx}\right)}
- \frac{\i}{n}B\left(
\zeta'(\bt).*\bx\right),
\]
then the value of the function $\mathcal{F}(\bx)$ can be computed quickly and accurately using the MATLAB function {\tt zfmm2dpart} in the MATLAB toolbox FMMLIB2D
developed by Greengard \& Gimbutas~\cite{Gre-Gim12}. This method for computing the solution $\psi$ to the BIE~\eqref{eq:ie-psi} is summarized in the following MATLAB function where the tolerances for the FMM and the GMRES method are taken to be $0.5\times10^{-15}$ and $10^{-14}$, respectively, and the GMRES method is run without restart:   
\begin{lstlisting}
	function y = szegofmm (et,etp,psi,n)
	Tet      =  etp./abs(etp);  Tet(etp==0) = 0;
	a        = [real(et.') ; imag(et.')];
	m        =  length(et)/n-1;
	y        =  gmres(@(x)F(x),psi,[],1e-14,100);
	function y = F(x)
	b1       = [abs(etp).*conj(x)].';
	[Ub1]    =  zfmm2dpart(5,(m+1)*n,a,b1,1);
	Eb1      = (Ub1.pot).';
	b2       = [abs(etp).*Tet.*x].';
	[Ub2]    =  zfmm2dpart(5,(m+1)*n,a,b2,1);
	Eb2      = (Ub2.pot).';
	y        =  x+(1./(n*i)).*(-conj(Tet).*conj(Eb1)+Eb2);
	end
	end
\end{lstlisting}

The preceding method assumes that the boundaries of the domains of interest are smooth, i.e. without corners. In the case of domains with corners (excluding cusps), to obtain accurate results, a re-parametrization  of the boundary of the domain is performed as described in~\cite{LSN17}.
Assume that the boundary component $\Gamma_j$ has $\ell$ corner points.
We first parametrize each boundary component $\Gamma_j$ by a $2\pi$-periodic function $\hat\zeta_j(t)$ for $t\in J_j=[0,2\pi]$. The function $\zeta_j(t)$ is assumed to be smooth with $\zeta'_j(t)\ne0$ for all values of $t\in J_j$ such that $\zeta_j(t)$ is not a corner point. 
We assume that $\zeta'_j(t)$ has only the first kind discontinuity at these corner points. At each corner point, the left tangent vector is taken to be the tangent vector at this point. As above, let $J$ be the disjoint union of the $m$ intervals $J_j=[0,2\pi]$, $j=1,2,\ldots,m$ and $\hat\zeta(t)$, $t\in J$, be a parametrization of the whole boundary $\Gamma$. 
Then, we parametrize the boundary $\Gamma$ by $\zeta(t)=\hat\zeta(\delta(t))$, where the function $\delta(t)$ is defined in~\cite[pp.~696--697]{LSN17}.
The function $\delta(t)$ is chosen such that the singularity in the first derivative of the solution of the BIE in the vicinity of the corner points is removed~\cite{Kre,LSN17}.
With the new parametrization $\zeta(t)$, the BIE~\eqref{eq:ie-psi} can be solved accurately using the above MATLAB function. 
However, for domains with corners, we usually need a larger number of points $n$ (which should be a multiple of the number of corners on each boundary component) for discretizing the BIE compared to domains with smooth boundaries (see~\cite{Kre,LSN17,Nas-ETNA} for further details).

Once the solution $\psi$ of the BIE~\eqref{eq:ie-psi} has been found, we can proceed to compute the analytic capacity $\gamma(E)$ using the formula~\eqref{eq:cap-mu2n}. This can be undertaken using the following MATLAB function:

\begin{lstlisting}
	function cap = ancap(zet,zetp,n)
	h    =  2*pi/n;
	rzet =  1/(2*pi)+zeros(size(zet));
	psi  =  szegofmm(zet,zetp,rzet,n);
	cap  =  sum(h*real(psi).*abs(zetp));
	end
\end{lstlisting}

\noindent
Various numerical examples will be presented in the proceeding two sections. We will take, in turn, domains bounded by Jordan curves and domains bounded by slits.

\section{Domains bounded by Jordan curves} \label{sec:Jordan}

In this section, we will use the method presented in the previous section to compute numerical approximations $\tilde\gamma(E)$ to the analytic capacity $\gamma(E)$ of compact sets bounded by smooth and piecewise smooth boundaries.

\begin{example}\label{ex:twodisks} 
	Consider the compact set $E=E_1\cup E_2$ where $E_{1,2}=\{z \in \CC\,|\,|z\pm c|\le r\}$ and $0<r<c$. Then~\cite{YR13}
	\[
	\gamma(E)=\frac{r}{2\sqrt{q}}\left(1-q\right)\theta_2^2(q), \quad q=\frac{p-\sqrt{p^2-1}}{p+\sqrt{p^2-1}}, \quad p=\frac{c}{r},
	\]
	where 
	\[
	\theta_2(q)=2q^{1/4}\prod_{j=1}^{\infty}(1-q^{2j})(1+q^{2j})^2
	\]
	is the second Jacobi theta function. 	
	The relative error in the computed approximate values $\tilde\gamma(E)$ obtained with $n=2^9$ are given in Table~\ref{tab:two-disks2} for several values of $c$ and $r$. 
	For $c=2$ and $r=1$, our obtained value is $\tilde\gamma(E)=1.875595019097120$ which is in the interval $(1.875595019097112,1.875595019097164)$ given in~\cite{YR13}. The exact value is
	\[
	\gamma(E)=\sqrt{3}\,\theta_2^2((2-\sqrt{3})^2)\approx 1.8755950190971197.
	\]
\end{example}

\begin{table}[h]
	\caption{The relative error in the approximate values of the analytic capacity $\gamma(E)$ for Example~\ref{ex:twodisks}.}
	\label{tab:two-disks2}
	\centering
	\begin{tabular}{l|cccc}  \hline
		$r$   & $c=0.5$  & $c=1$ & $c=2$ & $c=3$  \\ \hline
		$0.1$ & $1.40\times10^{-16}$     & $2.78\times10^{-16}$  & $1.39\times10^{-16}$    & $2.78\times10^{-16}$\\
		$0.5$ &       &$1.18\times10^{-16}$ &$1.80\times10^{-15}$ &$1.12\times10^{-15}$ \\
		$1$   &       &  &$1.18\times10^{-16}$ &$1.03\times10^{-15}$ \\
		$2$   &       &  &  &$2.99\times10^{-15}$ \\
		\hline
	\end{tabular}
\end{table}

\begin{example}\label{ex:square} 
	Let $E$ be the square with corners $1,-\i,-1,\i$. In this case, the domain $G = \widehat{\CC} \backslash E$ is an unbounded simply connected domain. Thus, by~\eqref{eq:cap-sq}, the analytic capacity of $E$ is 
	\begin{equation}\label{eq:cap-sq2}
		\gamma(E)=\capa(E)=\frac{\Gamma^2(1/4)}{2\sqrt{2\pi^3}}\approx 0.834626841674073.
	\end{equation}
	In this example, the boundary components of the compact set $E$ have corners. We use our numerical method with various values of $n$ to  approximate the analytic capacity $\gamma(E)$ and the results are presented in Table~\ref{tab:square}. As can clearly be seen in Table~\ref{tab:square}, the relative error decreases as $n$ increases. 
\end{example}

\begin{table}[h]
	\caption{The approximate values $\tilde\gamma(E)$ to the analytic capacity $\gamma(E)$ and their relative errors for Example~\ref{ex:square}, where $\gamma(E)$ is given by~\eqref{eq:cap-sq2}.}
	\label{tab:square}
	\centering
	\begin{tabular}{l|ccc}  \hline
		$n$              & $\tilde\gamma(E)$         & Relative error   & Time (sec)  \\ \hline
		$2^{8}$   & $0.834627510939279$   & $8.02\times10^{-7}$   & $0.13$   \\
		$2^{9}$   & $0.834626889228525$   & $5.70\times10^{-8}$   & $0.31$   \\
		$2^{10}$  & $0.834626845029913$   & $4.02\times10^{-9}$   & $0.40$   \\
		$2^{11}$  & $0.834626841908855$   & $2.81\times10^{-10}$  & $0.42$   \\
		$2^{12}$  & $0.834626841690432$   & $1.96\times10^{-11}$  & $0.55$   \\
		$2^{13}$  & $0.834626841675243$   & $1.40\times10^{-12}$  & $0.82$   \\
		$2^{14}$  & $0.834626841674219$   & $1.75\times10^{-13}$  & $1.02$   \\
		$2^{15}$  & $0.834626841674058$   & $1.80\times10^{-14}$  & $1.85$   \\
		\hline
	\end{tabular}
\end{table}

\begin{figure}[htb] %
	\centerline{
		\scalebox{0.5}{\includegraphics[trim=0 0 0 0,clip]{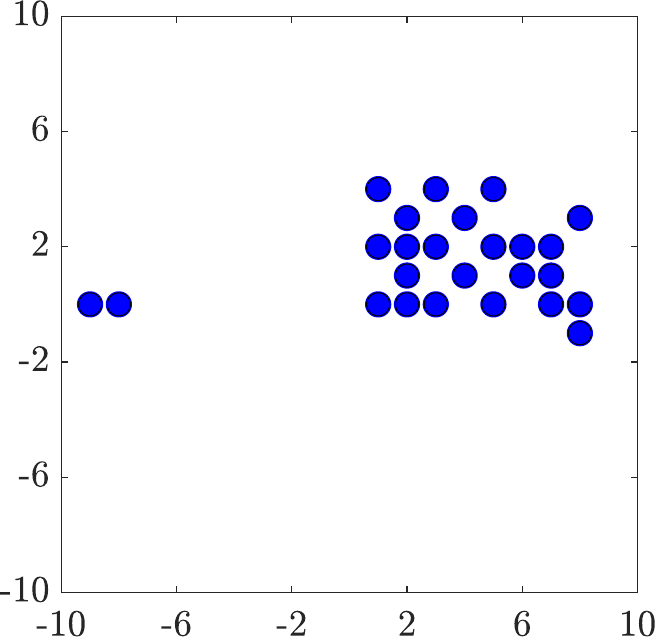}}
		\hfil
		\scalebox{0.5}{\includegraphics[trim=0 0 0 0,clip]{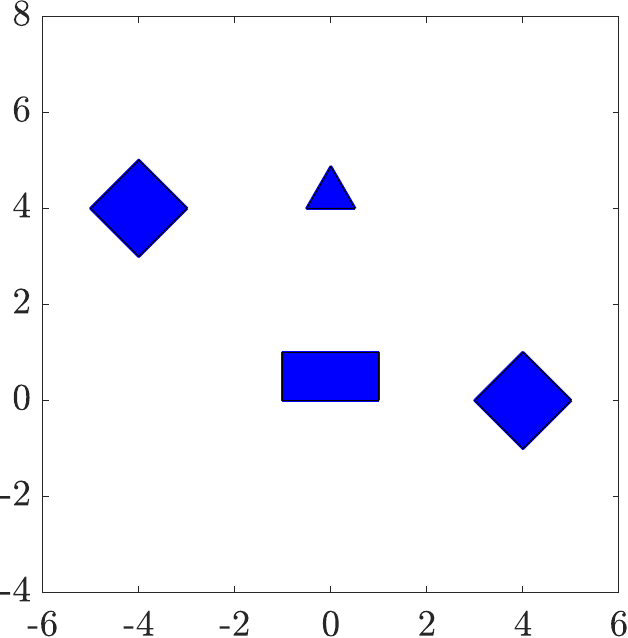}}
	}
	\centerline{
		\scalebox{0.5}{\includegraphics[trim=0 0 0 0,clip]{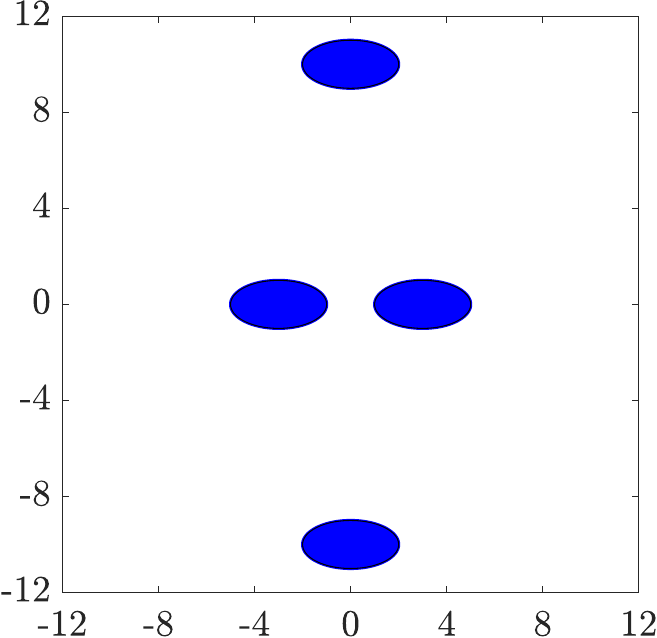}}
		\hfil
		\scalebox{0.5}{\includegraphics[trim=0 0 0 0,clip]{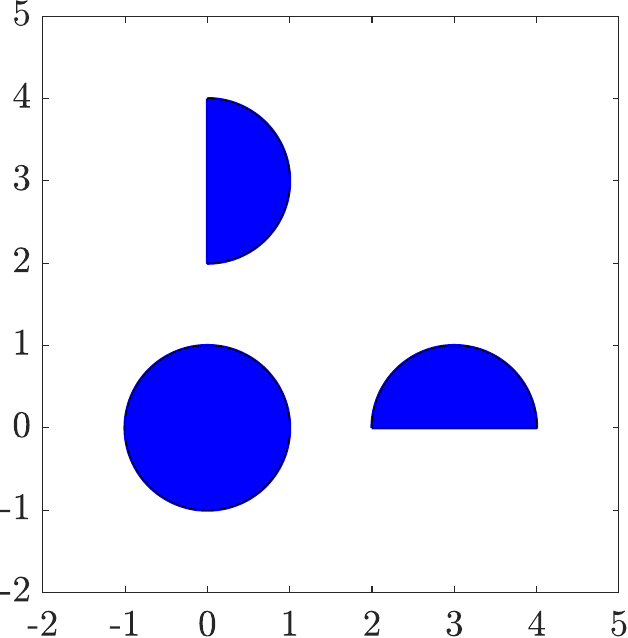}}
	}
	\caption{The four compact sets considered in Example~\ref{ex:4sets} and in Younsi \& Ransford~\cite[Figures~2, 3, 5 \& 6]{YR13}.}
	\label{fig:4sets}
\end{figure}

\begin{example}\label{ex:4sets} 
	As a validation of our numerical method, let us also consider the four compact sets shown in Figure~\ref{fig:4sets}. These sets were considered in~\cite[Figures~2, 3, 5 \& 6]{YR13}.
	We are not aware that the analytic capacities $\gamma(E)$ of these sets are known analytically. The approximate values of $\gamma(E)$ computed by our method are presented in Table~\ref{tab:4sets}. 
\end{example}

\begin{table}[h]
	\caption{The approximate values $\tilde\gamma(E)$ of the analytic capacity $\gamma(E)$ for Example~\ref{ex:4sets}.}
	\label{tab:4sets}
	\centering
	\begin{tabular}{llcc}  \hline
		$E$                     & $n$       & $\tilde\gamma(E)$          & Time (sec)  \\ \hline
		Union of $25$ disks     & $2^9$     & $4.14833193431691$   & $1.29$   \\
		Union of $4$ ellipses   & $2^9$     & $5.37199577528044$   & $0.55$   \\		
		Union of $4$ polygons   & $3\times2^{12}$  & $2.69399397476757$   & $3.28$   \\		
		Union of a disk \& $2$ semi-disks & $2^{10}$  & $2.12070613733366$   & $0.57$   \\		
		\hline
	\end{tabular}
\end{table}

It is important to point out that for the union of $4$ ellipses, the elapsed time of $0.55$ seconds suggests that our computations are more than $1000$ times faster than the method used for the same problem in~\cite{YR18,YR13}. However, around a decade has passed since the computations in~\cite{YR13} were performed. We also note that our method is much faster than the method presented in~\cite{YR18,YR13} for non-circular compact sets compared with circular ones.

\begin{figure}[htb] %
	\centerline{
		\scalebox{0.375}{\includegraphics[trim=1.8cm 0 1.8cm 0,clip]{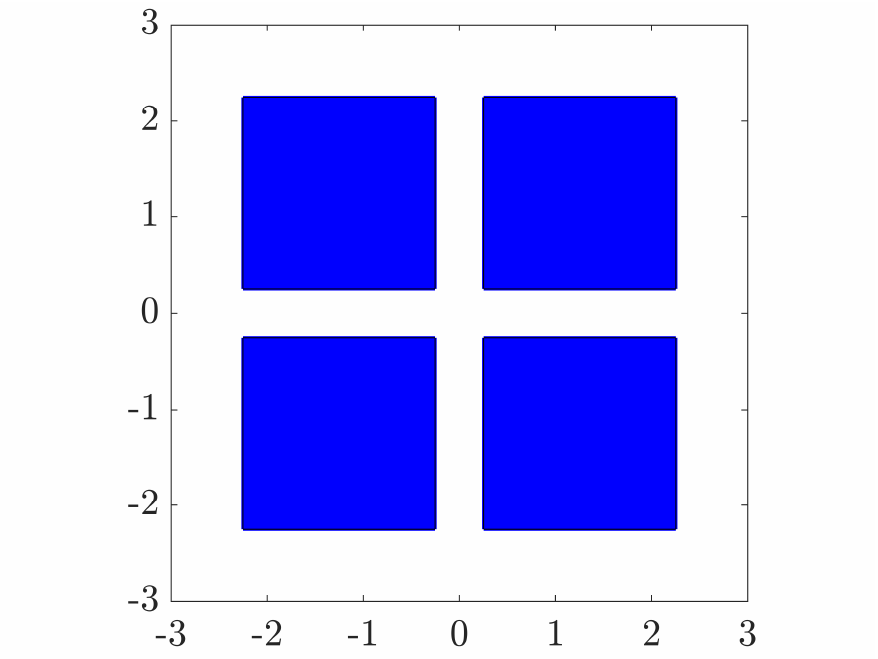}}
		\hfil
		\scalebox{0.375}{\includegraphics[trim=1.8cm 0 1.8cm 0,clip]{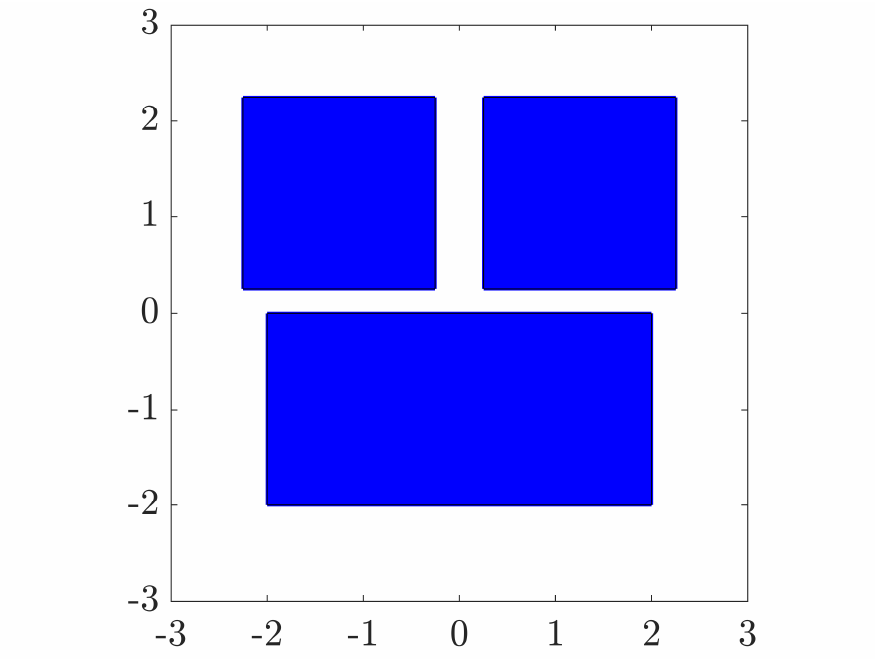}}
		\hfil
		\scalebox{0.375}{\includegraphics[trim=1.8cm 0 1.8cm 0,clip]{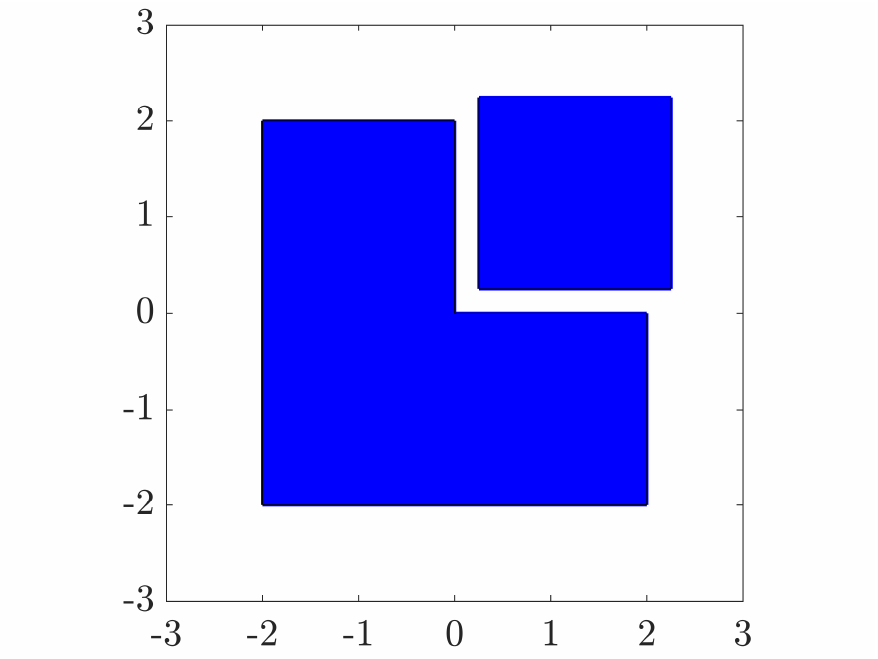}}
	}
	\caption{The three kinds of compact sets considered in Example~\ref{ex:4squares}.}
	\label{fig:4squares}
\end{figure}

\begin{figure}[htb] %
	\centerline{
		\scalebox{0.5}{\includegraphics[trim=0 0 0 0,clip]{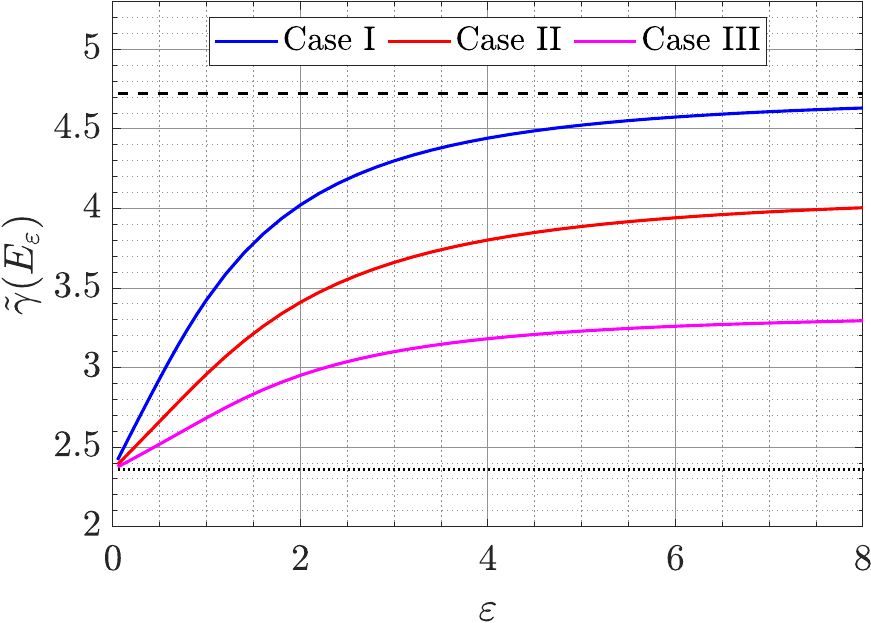}}
	}
	\caption{Graphs of the approximate values $\tilde\gamma(E_\varepsilon)$ of the analytic capacity as a function of $\varepsilon$ in each of the three cases in Example~\ref{ex:4squares}. The dashed line denotes the upper bound $2\gamma(E_0)$ and the dotted line denotes the lower bound $\gamma(E_0)$.}
	\label{fig:4sqcap}
\end{figure}

\begin{example}\label{ex:4squares} 
	Consider the square $[-2,2]\times[-2,2]$ and the four sub-squares with centers $\pm 1 \pm \i$. 
	Let $\varepsilon>0$ be a real parameter.  We consider three cases (i)-(iii). In the case (i), consider moving the centers of the sub-squares via the parameter $\epsilon$ to the points $(1+\varepsilon)(\pm 1 \pm \i)$ (see Figure~\ref{fig:4squares} (left)). In case (ii), let us fix the lower two sub-squares and consider moving the centers of the upper sub-squares to the points $(1+\varepsilon)(\pm 1 + \i)$ (see Figure~\ref{fig:4squares} (middle)). Finally, in case (iii), let us fix three of the sub-squares and move the center of the remaining sub-square to the point $(1+\varepsilon)(1+\i)$ by increasing the value of $\epsilon$ (see Figure~\ref{fig:4squares} (right)).  
	Let us label the union of the compact sets generated with $\varepsilon$ in each case by $E_\varepsilon$. For $\varepsilon=0$, then $E_0=[-2,2]\times[-2,2]$ is the original square. 
	
	Note that the length of each side of the original square is $4$, and hence, by~\eqref{eq:cap-sq},
	\[
	\gamma(E_0)=\frac{\Gamma^2(1/4)}{\sqrt{\pi^3}}.
	\]
	Further, let us label the sub-squares by $F_1,\ldots,F_4$, then the length of each of these sub-squares is $2$, and hence, by~\eqref{eq:cap-sq},
	\[
	\gamma(F_j)=\frac{\Gamma^2(1/4)}{2\sqrt{\pi^3}}=\frac{1}{2}\gamma(E_0), \quad j=1,\ldots,4.
	\]
	In each of the three cases, we compute the analytic capacity $\gamma(E_\varepsilon)$ as a function of $\varepsilon$ which are presented in Figure~\ref{fig:4sqcap}.
	The results indicate that $\gamma(E_0)$ is a lower bound for $\gamma(E_\varepsilon)$ and $\sum_{j=1}^4\gamma(F_j)=2\gamma(E_0)$ is an upper bound for $\gamma(E_\varepsilon)$. These results collectively provide numerical evidence to corroborate the conjectural subadditivity property of analytic capacity for these compact sets. Our results demonstrate the expected phenomenon that the analytic capacity increases as the sub-squares move further apart from each other, and values of the analytic capacity in the case (i) are larger than those in the cases (ii) and (iii) corresponding to the greater number of sub-squares.
	
\end{example}

\begin{figure}[htb] %
	\centerline{
		\scalebox{0.375}{\includegraphics[trim=1.8cm 0 1.8cm 0,clip]{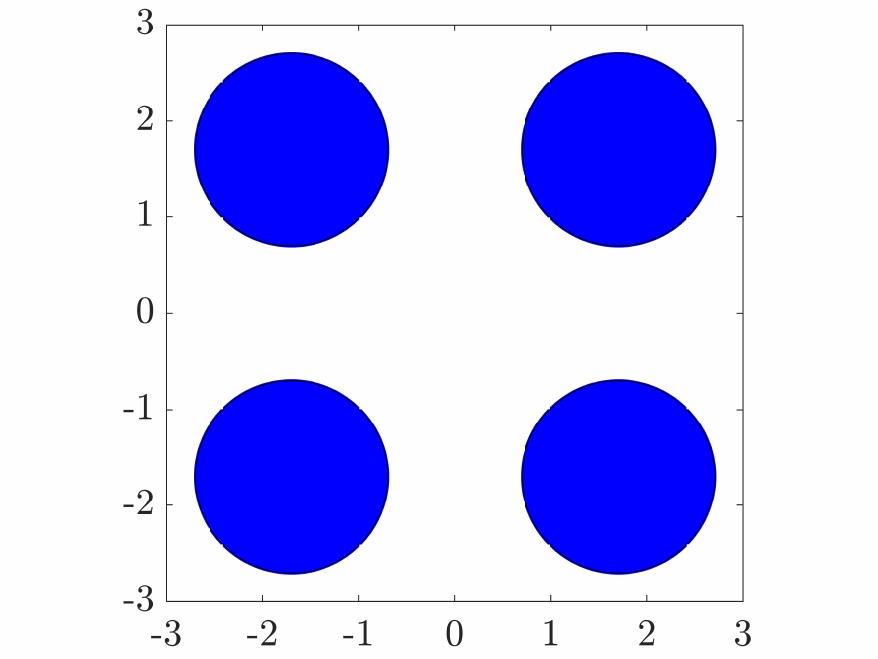}}
		\hfil
		\scalebox{0.375}{\includegraphics[trim=1.8cm 0 1.8cm 0,clip]{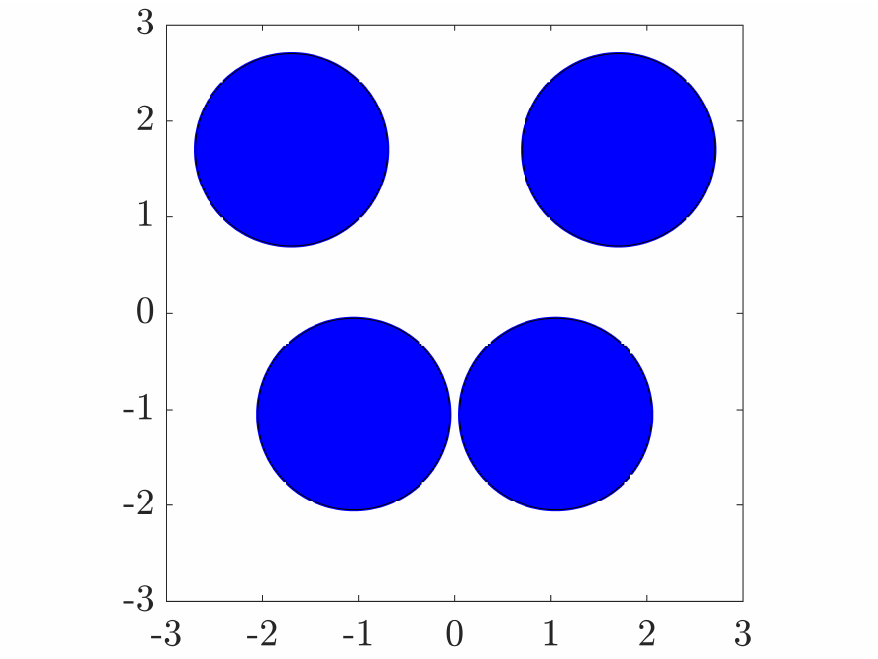}}
		\hfil
		\scalebox{0.375}{\includegraphics[trim=1.8cm 0 1.8cm 0,clip]{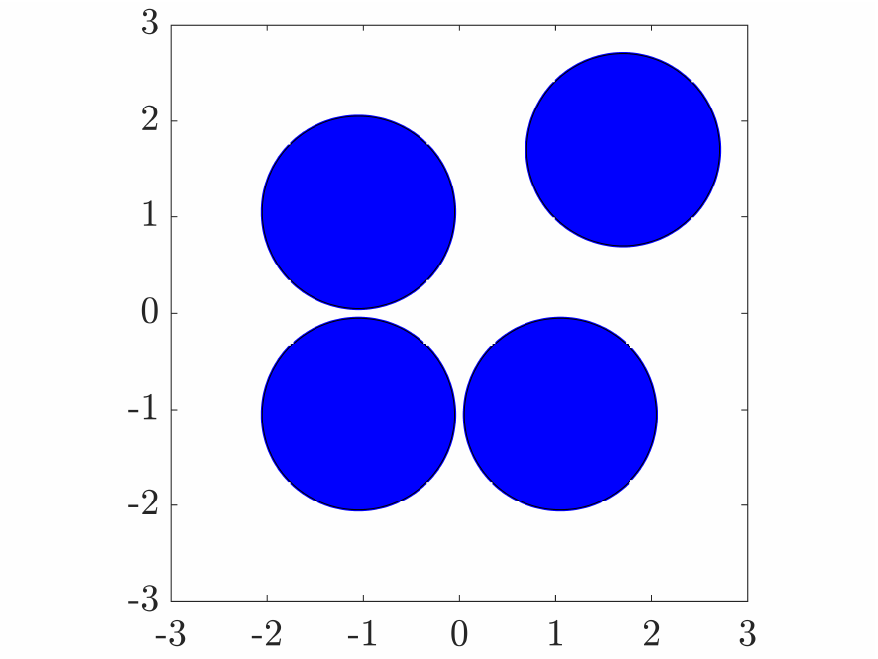}}
	}
	\caption{The three kinds of compact sets considered in Example~\ref{ex:4squares}.}
	\label{fig:4disks}
\end{figure}

\begin{figure}[htb] %
	\centerline{
		\scalebox{0.5}{\includegraphics[trim=0 0 0 0,clip]{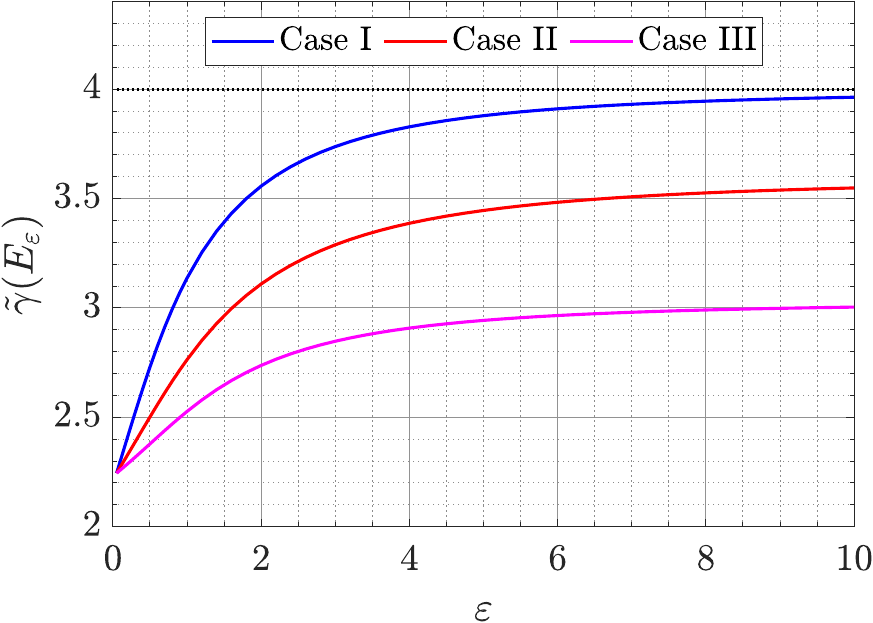}}
	}
	\caption{Graphs of the approximate values $\tilde\gamma(E_\varepsilon)$ of the analytic capacity as a function of $\varepsilon$ in each of the three cases in Example~\ref{ex:4disks}. The dotted line denotes the upper bound $\sum_{j=1}^4\gamma(F_j)=4$.}
	\label{fig:4diskcap}
\end{figure}

\begin{example}\label{ex:4disks} 
	Consider the four disks with radius $1$ and centers at $\pm1\pm\i$. Let $\varepsilon>0$ be a real parameter.  We consider three cases (i)-(iii). In case (i), consider moving the centers of the four disks via the parameter $\epsilon$ to the points $(1+\varepsilon)(\pm 1 \pm \i)$ (see Figure~\ref{fig:4disks} (left)). In case (ii), let us fix the lower two disks and consider moving the centers of the upper disks to the points $(1+\varepsilon)(\pm 1 + \i)$ (see Figure~\ref{fig:4disks} (middle)). Finally, in case (iii), let us fix three of the disks and move the center of the remaining disk to the point $(1+\varepsilon)(1+\i)$ by increasing the value of $\epsilon$ (see Figure~\ref{fig:4disks} (right)).  
	We label the union of the compact sets generated with $\varepsilon$ in each case by $E_\varepsilon$. 
	
	Let us label the disks by $F_1,\ldots,F_4$, then $\gamma(F_j)=1$, $j=1,\ldots,4$. In each of the three cases, we compute the analytic capacity $\gamma(E_\varepsilon)$ as a function of $\varepsilon$ which are presented in Figure~\ref{fig:4diskcap}. The results demonstrate the expected phenomenon that the analytic capacity increases as the disks move further apart from each other, and values of the analytic capacity in case (i), corresponding to the case when all disks are moving away from each other, are larger than those in cases (ii) and (iii).	
\end{example}

\begin{remark}\label{rem:disks}
	If $E$ is a compact set such that $G = \widehat{\CC} \backslash E$ is bounded by finitely many analytic curves, then~\cite{YR18,YR13}
	\[
	\gamma(E)=\min\left\{\frac{1}{2\pi}\int_{\partial G}|1+g(z)|^2|dz|\right\}
	\] 	
	where the minimum is taken over all functions $g(z)$ that are analytic on $G$, continuous up to the boundary, with $g(\infty)=0$. Thus, taking $g(z)\equiv0$ gives the inequality
	\[
	\gamma(E)\le \frac{1}{2\pi}\int_{\partial G}|dz|.
	\]
	Hence, if $E$ is a union of finitely many disjoint closed disks $E_1,\ldots,E_m$ with radii $r_1,\ldots,r_m$, we have
	\[
	\gamma(E)\le r_1+\cdots+r_m=\gamma(E_1)+\cdots+\gamma(E_m).
	\]
	This explains the results presented in Figure~\ref{fig:4diskcap} indicating that $\sum_{j=1}^4\gamma(F_j)=4$ is an upper bound for $\gamma(E_\varepsilon)$. 
\end{remark}

\begin{figure}
	\centerline{
		\scalebox{0.45}{\includegraphics[trim=0 0 0 0,clip]{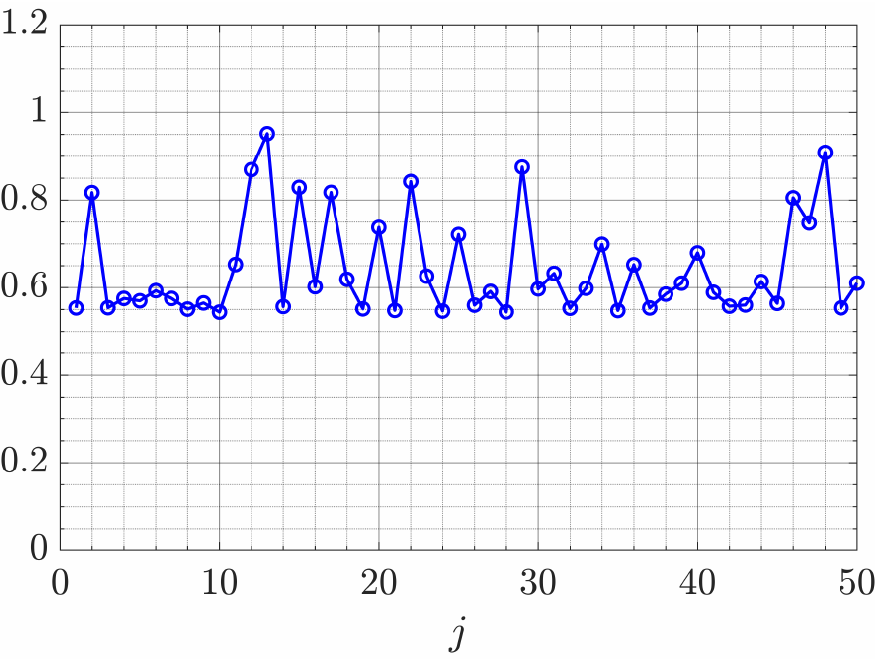}}
		\hfil
		\scalebox{0.41}{\includegraphics[trim=1cm -1cm 0 0,clip]{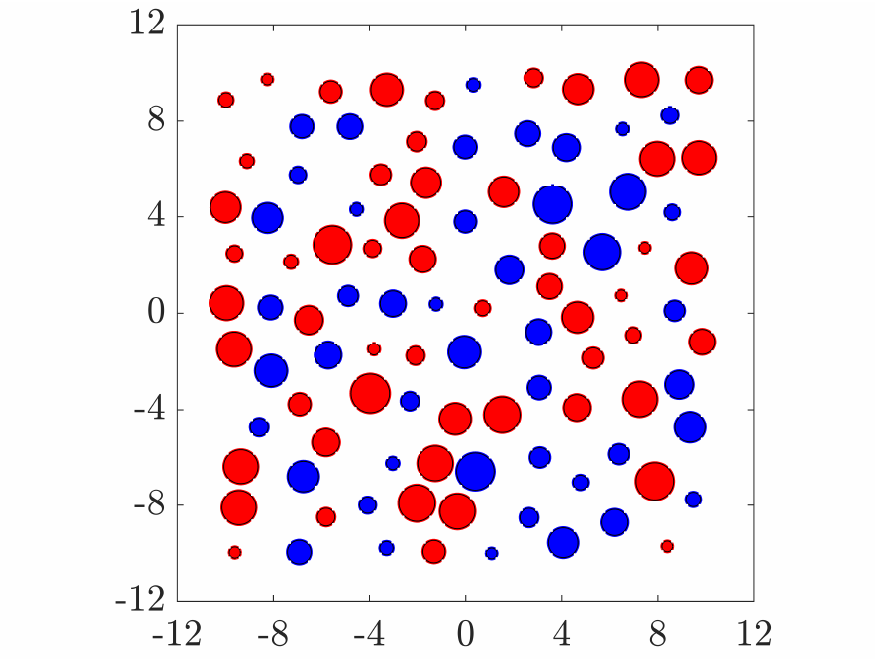}}
	}
	\caption{The ratio in~\eqref{eq:ratios} as functions of the run number $j$ for random radii in $(0.2,0.8)$ (left), and an example of the compact set $E \cup F$ with $100$ random disks, where the elements of $E$ are shown in blue and the elements of $F$ are shown in red (right).}
	\label{fig:100disks}
\end{figure}

\begin{example}\label{ex:100disks} 
	We consider $m=100$ random non-overlapping disks. For $k=1,2,\ldots,100$, the radius $r_k$ of the disk $E_k$ is chosen randomly in $(0.2,0.8)$ and its center $c_k$ is chosen in the square  $[-10,10]\times[-10,10]$ such that all disks are non-overlapping. Then we randomly choose an integer $\ell\in[1,99]$. We define 
	\[
	E=\bigcup_{k=1}^{\ell}E_k,\quad F=\bigcup_{k=\ell+1}^{m}E_k.
	\]  
	See Figure~\ref{fig:100disks} (right) for an example of such compact sets $E$ and $F$. 
	
	For this problem, we run our method $50$ times to obtain $50$ different locations for these disks as well as different sets $E$ and $F$. For each run $j$, we use the above presented method with $n=2^9$ to compute approximate values of the quantities $\gamma(E)$, $\gamma(F)$ and $\gamma(E\cup F)$. The values of the ratio 
	\begin{equation}\label{eq:ratios}
		\frac{\gamma(E\cup F)}{\gamma(E)+\gamma(F)}
	\end{equation}
	are plotted as a function of the run number $j$ as shown in Figure~\ref{fig:100disks} (left). 
	As can be seen from the graphs in this figure, we have verified that the conjectural subadditivity property of analytic capacity holds for each of the 50 random compact sets we considered, and in particular that
	\[
	\gamma(E\cup F)\le \gamma(E)+\gamma(F).
	\]  
	Again, the conjectural subadditivity property for analytic capacity holds for the compact sets in this example. 
	
\end{example}

\begin{figure}[htb] %
	\centerline{
		\scalebox{0.45}{\includegraphics[trim=0 0 0 0,clip]{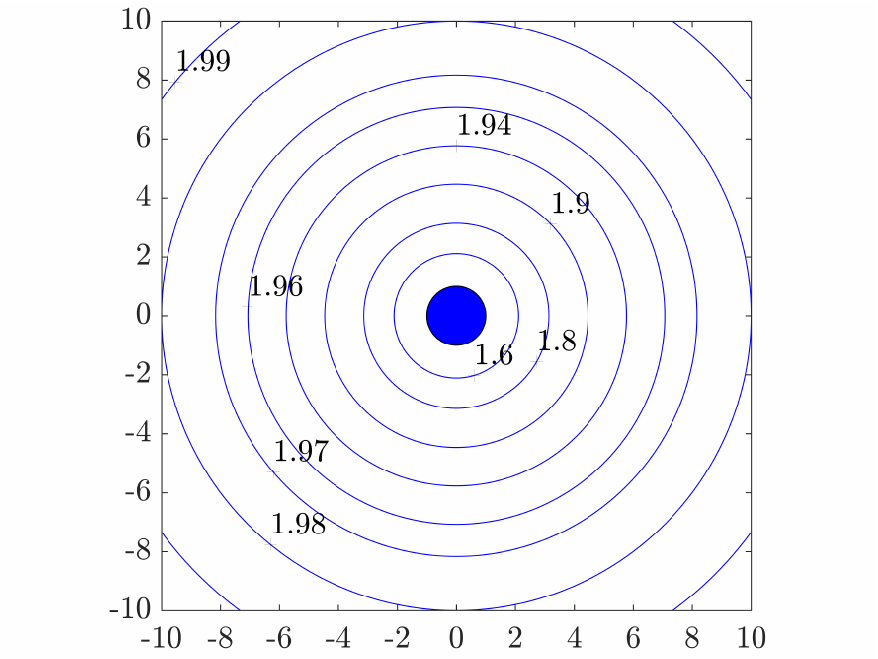}}
		\hfil
		\scalebox{0.45}{\includegraphics[trim=0 0cm 0 0,clip]{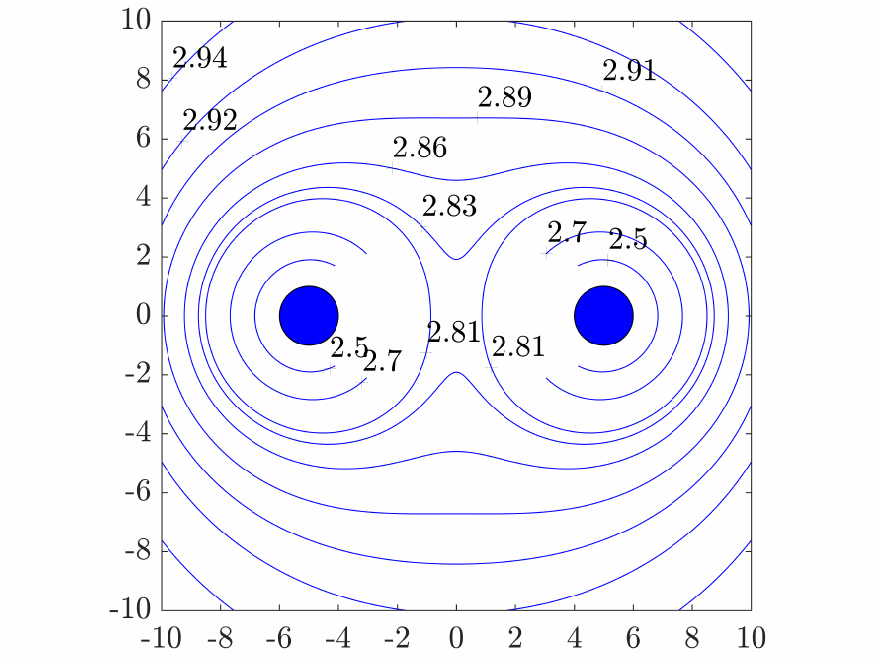}}
	}
	\centerline{
		\scalebox{0.45}{\includegraphics[trim=0 0 0 0,clip]{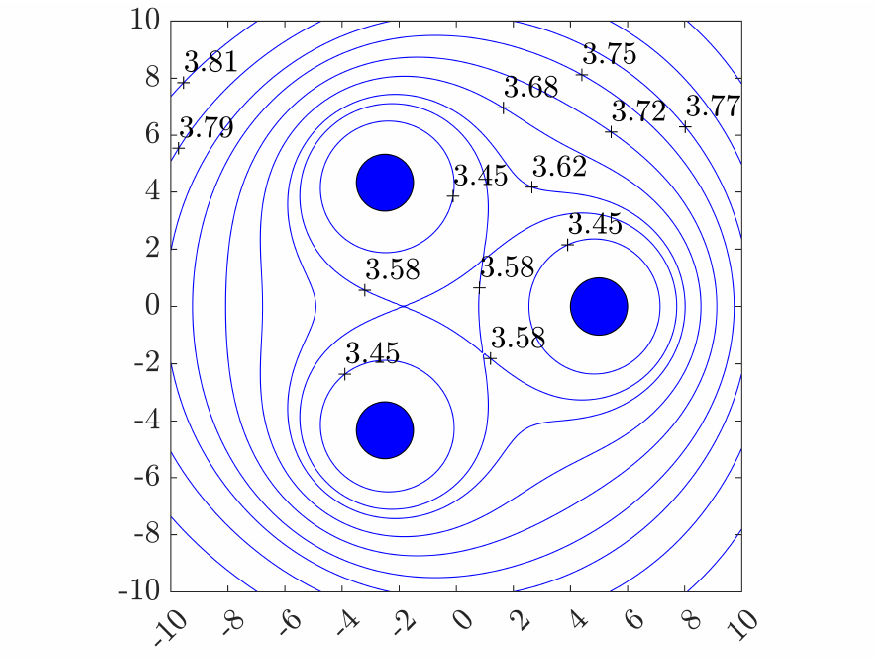}}
		\hfil
		\scalebox{0.45}{\includegraphics[trim=0 0cm 0 0,clip]{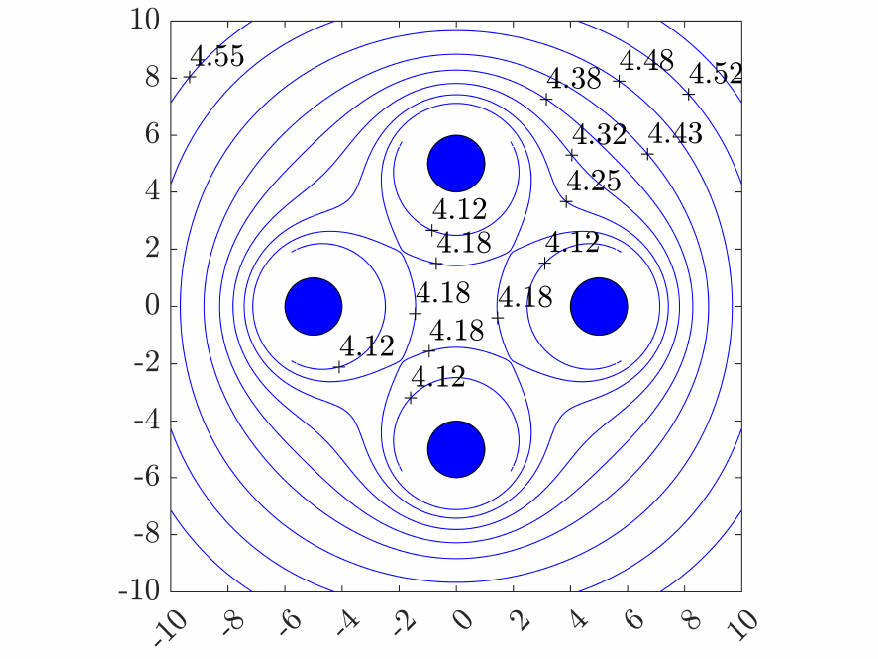}}
	}
	\caption{The level curves of the function $u(x,y)=\gamma(E_m\cup F)$ in Example~\ref{ex:mdisks} when $m=1,2,4,5$, and the disks forming $E_m$ are colored blue.}
	\label{fig:level}
\end{figure}

\begin{figure}
	\centerline{
		\scalebox{0.45}{\includegraphics[trim=0 0 0 0,clip]{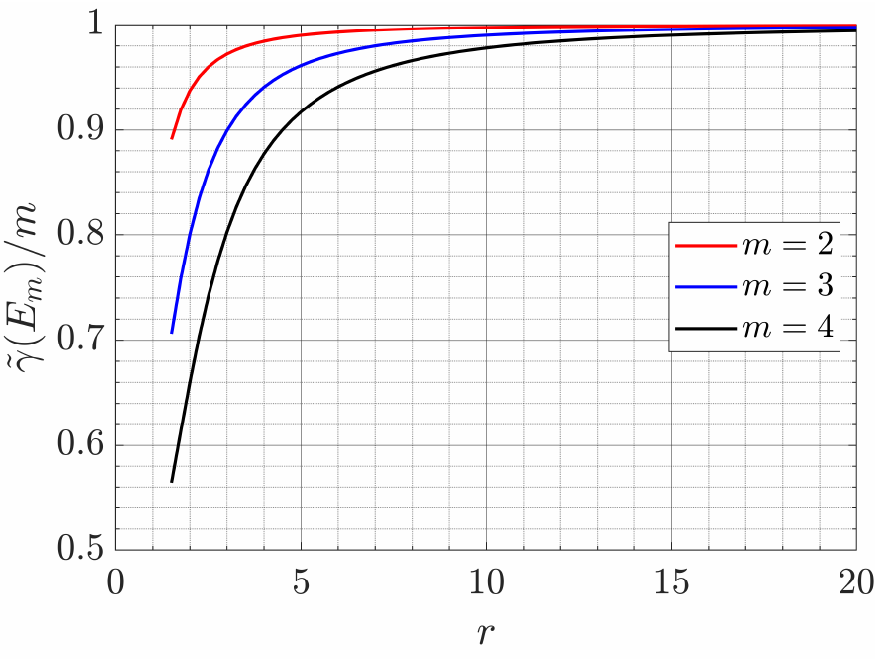}}
	}
	\caption{The ratios $\tilde\gamma(E_m)/m$ as functions of the distance $r$ measured from the origin to the centers of the $m$ disks.}
	\label{fig:4disksr}
\end{figure}

\begin{example}\label{ex:mdisks} 
	In this example, we study compact sets consisting of disjoint disks of equal radii. In particular, we will validate numerically the conjectural subadditivity property of analytic capacity in several cases. The consideration of such compact sets is important since proving the conjectural subadditivity property of analytic capacity for arbitrary compact sets $E,F\in\CC$ is equivalent to proving it for all disjoint compact sets that are finite unions of disjoint closed disks, all with the same radius~\cite{Mel,YR13}.

	Let $E_m$ be a union of $m$ disjoint disks and $F$ be a disk with center $x+\i y$ such that these $m+1$ disks are non-overlapping and of unit radii. Let the real function $u(x,y)$ be defined by
	\[
	u(x,y) = \gamma(E_m\cup F),
	\] 	
	i.e. the function $u(x,y)$ is defined for all points $(x,y)\in\RR^2$ such that the disk $F$ does not overlap any of the other $m$ disks. Note that the domain of definition of the function $u(x,y)$ is unbounded. Hence, in our numerical computations, we consider only the region $-10\le x,y\le 10$ and we assume that the minimum distance between any two disks is $0.02$. We compute approximate values of the function $u(x,y)$ using $n=2^{10}$ and then plot several level curves of $u(x,y)$. 
	We consider four cases of $m$, namely $m=1,2,3,4$. For $m=1$, $E_m$ consists of only one disk which is assumed to be the unit disk. For $m=2,3,4$, we assume that the centers of the disks forming $E_m$ are $5e^{2k\pi\i/m}$, $k=1,2,\ldots,m$. The approximate values $\tilde\gamma(E_m)$ of the analytic capacity of $E_m$ are presented in the following table.

	\begin{center}
		\begin{tabular}{l@{\qquad}l}  \hline
			$m$     & $\tilde\gamma(E_m)$    \\ \hline
			$1$     & $1$      \\
			$2$     & $1.98000206142844$      \\
			$3$     & $2.88420404308815$      \\
			$4$     & $3.67012955644439$      \\
			\hline
		\end{tabular}
	\end{center}
	
	\noindent
	The computed level curves of the function $u(x,y)$  are presented in Figure~\ref{fig:level}. 
	The results presented demonstrate  that the values of the analytic capacity $\gamma(E_m\cup F)$ increase when $F$ moves away from $E_m$ and decrease when $F$ moves towards $E_m$. In particular, when $F$ is far away from $E_m$, we have 
	\begin{equation}\label{eq:Em-F}
		\gamma(E_m\cup F) \approx \gamma(E_m)+\gamma(F).
	\end{equation}
	The results presented also illustrate that 
	\[
	\gamma(E_m) \le \gamma(E_m\cup F) \le \gamma(E_m)+\gamma(F)=\gamma(E_m)+1,
	\]
	which indicates that the conjectural subadditivity property for analytic capacity again holds in this example.  	
	
	Finally, we compute approximate values $\tilde\gamma(E_m)$ of the analytic capacity for the above compact sets $E_m$, this time assuming that the centers of the disks are $re^{2k\pi\i/m}$, $k=1,2,\ldots,m$, for $1.5\le r\le 20$. The ratios $\tilde\gamma(E_m)/m$ are plotted in Figure~\ref{fig:4disksr} as functions of $r$. 
	As we can see, each of these ratios tends to $1$ as $r$ increases; that is, the analytic capacity of $E_m$ (the union of $m$ disks) tends to the sum of the individual analytic capacities of these $m$ disks which is equal to $m$ (see Remark~\ref{rem:disks} above). The results presented in~\eqref{eq:Em-F} and Figures~\ref{fig:level} and~\ref{fig:4disksr} provide experimental evidence for a previously known result due to Pommerenke who showed that, roughly speaking, the analytic capacity of a compact set is approximately equal to the sum of the analytic capacities of its components if they are far away from each other (see also~\cite[p.~267]{Kir}).
	
\end{example}	

\section{Domains bounded by slits} \label{sec:Slits}

In this section, we will use our numerical method to compute numerical approximations $\tilde\gamma(E)$ to the analytic capacity $\gamma(E)$ of sets consisting of slits. We will consider here only rectilinear slits. However, the method presented can be extended to other types of slits.

For $a_1,\ldots,a_m,b_1,\ldots,b_m\in\CC$ such that $E_j=[a_j,b_j]$ are non-overlapping for $j=1,2,\ldots,m$, let $E=\bigcup_{j=1}^mE_j$ and let $\Omega=\widehat{\CC} \backslash E$, i.e., $\Omega$ is the unbounded multiply connected slit domain obtained by removing the $m$ slits $E_1,\ldots,E_m$ from the extended complex plane $\overline{\CC}$. The method presented in Section~\ref{sec:nm} is \emph{not} directly applicable to such a domain $\Omega$. 
However, an iterative method has been presented in~\cite{NG18} for constructing of a preimage unbounded multiply connected domain $G$ bounded by smooth Jordan curves and the unique conformal mapping $\zeta=\Phi(z)$, $\Phi\,:\,G\to\Omega$, such that $\Phi$ is normalized near infinity by the condition
\begin{equation}\label{eq:cm-recs}
	\Phi(z)=z+O\left(\frac{1}{z}\right).
\end{equation} 
The inverse function $z=\Phi^{-1}(\zeta)$ is then the conformal mapping from $\Omega$ onto $G$. The method presented in Section~\ref{sec:nm} is now applicable to the new domain $G$. This iterative method is reviewed in Appendix~\ref{sec:App-pre}. 

Let $w=f(z)$ be the Ahlfors map from the unbounded domain $G$ onto the unit disk $\DD$ such that $f(\infty)=0$, $f'(\infty)>0$, and $f'(\infty)$ is maximal. Then the function 
\[
w=g(\zeta)=f\left(\Phi^{-1}(\zeta)\right)
\]
is an Ahlfors map from the unbounded slit domain $\Omega$ onto the unit disk $\DD$. We have
\[
g(\infty)=f\left(\Phi^{-1}(\infty)\right)=f(\infty)=0,
\]
and
\[
g'(\infty)=\lim_{\zeta\to\infty}\zeta[g(\zeta)-g(\infty)]=\lim_{\zeta\to\infty}\zeta f\left(\Phi^{-1}(\zeta)\right).
\]
If $z=\Phi^{-1}(\zeta)$, then $\zeta\to\infty$ as $z\to\infty$, and 
\[
g'(\infty)=\lim_{z\to\infty}\Phi(z) f(z)=\lim_{z\to\infty}\frac{\Phi(z)}{z}\,\lim_{z\to\infty}zf(z)=f'(\infty),
\]
where we used~\eqref{eq:cm-recs} and the fact that $\lim_{z\to\infty}zf(z)=f'(\infty)$.
Since $f'(\infty)$ is maximal, it then follows that $g'(\infty)$ is maximal too.
Thus, by computing the preimage domain $G$ and the conformal mapping $\Phi\,:\,G\to\Omega$ normalized by the condition~\eqref{eq:cm-recs}, we will have
\[
\gamma(E)=g'(\infty)=f'(\infty)
\]
where $w=f(z)$ is the Ahlfors map from the unbounded domain $G$ onto the unit disk $\DD$ with the normalization $f(\infty)=0$ and $f'(\infty)>0$. Since the domain $G$ is bounded by smooth Jordan curves, the value of $f'(\infty)$ can be computed as explained in Section~\ref{sec:nm}. 

We consider three examples. In our first example, we consider rectilinear slits on the real line. We know the exact value of the analytic capacity in this case and hence the error in the computed approximate values can be calculated. Examples with unknown explicit formulae are also presented.

\begin{example}\label{ex:cant}
	We consider several of the compact sets used in the process of generating the middle-thirds Cantor set. Let $E_0=[-1,1]$, and let $E_k$ for $k\ge1$ be defined recursively by
	\begin{equation}\label{eq:Ek}
		E_k=\frac{1}{3}\left(E_{k-1}-\frac{1}{3}\right)\bigcup\frac{1}{3}\left(E_{k-1}+\frac{1}{3}\right).
	\end{equation}
	Note that $E_k$ consists of $m=2^k$ sub-intervals of the $[0,1]$, each of length $2/3^k$. We denote these sub-intervals by $I_j$ for $j=1,2,\ldots,m$. 
	By~\eqref{eq:cap-int}, the exact value of $\gamma(E_k)$ is known and given by
	\begin{equation}\label{eq:Ek-ex}
		\gamma(E_k)=\frac{1}{4}|E_k|=\frac{1}{4}\sum_{j=1}^{m}|I_j|=\frac{1}{2}\,\left(\frac{2}{3}\right)^k,
	\end{equation}
	from which it is immediate that $\gamma(E_k)\to0$ as $k\to\infty$. Note for this example that
	\[
	\gamma(E_k)=\gamma\left(\bigcup_{j=1}^{m}I_j\right)=\sum_{j=1}^{2^k}\gamma(I_j).
	\]
	The proposed method is used to compute approximate values $\tilde\gamma(E_k)$ to the analytic capacity $\gamma(E_k)$ for $k=1,2,\ldots,10$ and the obtained results are presented in Table~\ref{tab:RecSlit}. We compute also the relative error in the computed approximate values. 
	As can be seen in Table~\ref{tab:RecSlit}, our numerical method gives accurate results for sets consisting of a very high number of slits. 
\end{example}

\begin{table}[h]
	\caption{The approximate values $\tilde\gamma(E_k)$ to the analytic capacity $\gamma(E_k)$ and their relative errors for Example~\ref{ex:cant}, where $\gamma(E_k)$ is given by~\eqref{eq:Ek-ex}.}
	\label{tab:RecSlit}
	\centering
	\begin{tabular}{llll}
		\hline
		$k$   & $m$   & $\tilde\gamma(E_k)$   & Relative error  \\ \hline
		$1$   & $2$   & $0.333333333333328$   & $1.45\times10^{-14}$  \\
		$2$   & $4$   & $0.222222222222224$   & $7.62\times10^{-15}$  \\
		$3$   & $8$   & $0.148148148148148$   & $2.06\times10^{-15}$  \\
		$4$   & $16$  & $0.098765432098761$   & $4.29\times10^{-14}$  \\
		$5$   & $32$  & $0.065843621399160$   & $2.60\times10^{-13}$  \\
		$6$   & $64$  & $0.043895747599439$   & $2.74\times10^{-13}$  \\
		$7$   & $128$ & $0.029263831732889$   & $2.68\times10^{-12}$  \\
		$8$   & $256$ & $0.019509221155260$   & $2.66\times10^{-12}$  \\
		$9$   & $512$ & $0.013006147436531$   & $2.64\times10^{-11}$  \\
		$10$  & $1024$ & $0.008670764957687$   & $2.65\times10^{-11}$  \\
		\hline
	\end{tabular}
\end{table}

\begin{figure}[htb] %
	\centerline{
		\scalebox{0.5}{\includegraphics[trim=0 0 0 0,clip]{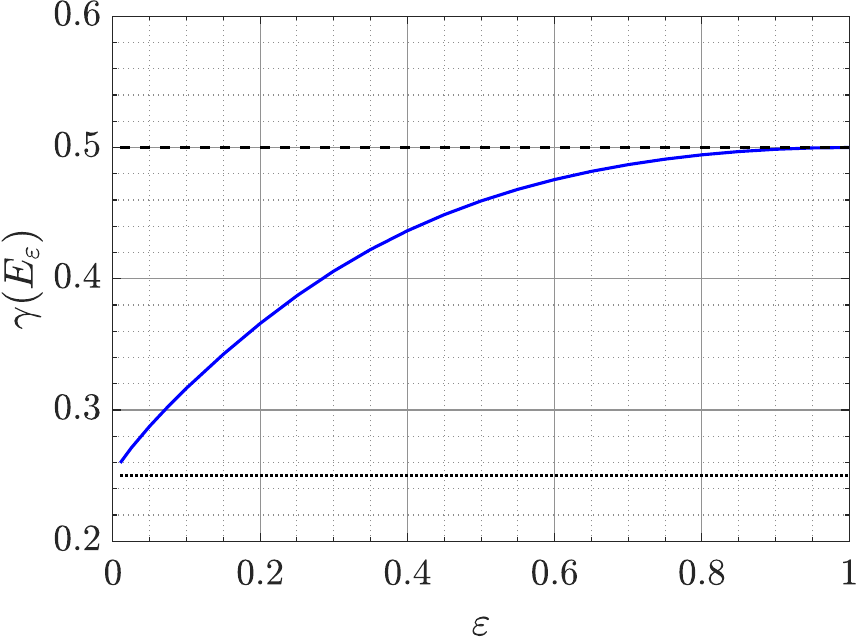}}
		\hfil
		\scalebox{0.455}{\includegraphics[trim=0 -1cm 0 0,clip]{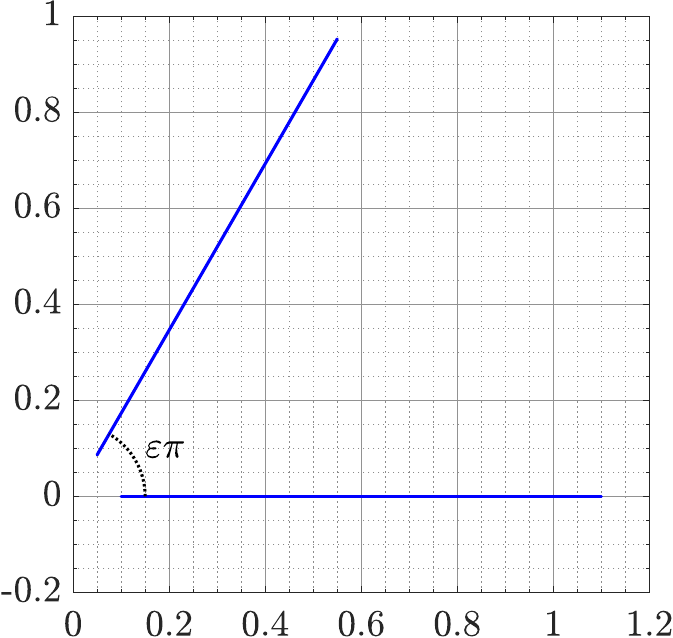}}
	}
	\caption{Graph of the approximate values $\tilde \gamma(E_{\varepsilon})$ of the analytic capacity for the compact set $E_\varepsilon$ in Example~\ref{ex:4sets} as a function of $\varepsilon$ (left). A schematic of the compact set $E_\varepsilon$ showing the parameter $\varepsilon$ (right).}
	\label{fig:flap}
\end{figure}

\begin{example}\label{ex:flap}
	We next consider the union of two equal rectilinear slits of unit length: one slit $F_1=[0.1,1.1]$ is fixed on the real line, and the other is taken to be $F_{2,\varepsilon}=[0.1e^{\i \varepsilon\pi},1.1e^{\i \varepsilon\pi}]$ where we vary $\varepsilon$ between zero and one. A schematic of this configuration is shown in Figure~\ref{fig:flap} (right) when $\varepsilon=1/3$. 	
	Let $E_\varepsilon=F_1 \cup F_{2,\varepsilon}$. It is known that $\gamma(F_1)=\gamma(F_{2,\varepsilon})=1/4$, by~\eqref{eq:cap-seg}.
	When $\varepsilon=1$, we have $E_1=[0.1,1.1]\cup[-1.1,-0.1]$ and hence $\gamma(E_1)=1/2$, by~\eqref{eq:cap-int}.
	When $\varepsilon=0$, we have $E_0=F_1=F_{2,0}=[0.1,1.1]$ and hence $\gamma(E_0)=1/4$.
	For $0<\varepsilon<1$, there is no exact value of $\gamma(E_\varepsilon)$. We use our method to compute $\gamma(E_\varepsilon)$ for $0.05\le\varepsilon\le1$ and the numerical results are presented in Figure~\ref{fig:flap} (left). It is clear that 
	\[
	1/4=\gamma(E_0)\le \gamma(E_\varepsilon) \le \gamma(E_1)=1/2.
	\]
	That is, the value of the analytic capacity of $E_\varepsilon$ is maximum when the two slits are collinear. Furthermore, we always have 
	\[
	\gamma(E_\varepsilon)=\gamma(F_1 \cup F_{2,\varepsilon}) \le \gamma(E_1)+\gamma(F_{2,\varepsilon})=1/2.
	\]
\end{example}

\begin{figure}[htb] %
	\centerline{
		\scalebox{0.5}{\includegraphics[trim=0 0 0 0,clip]{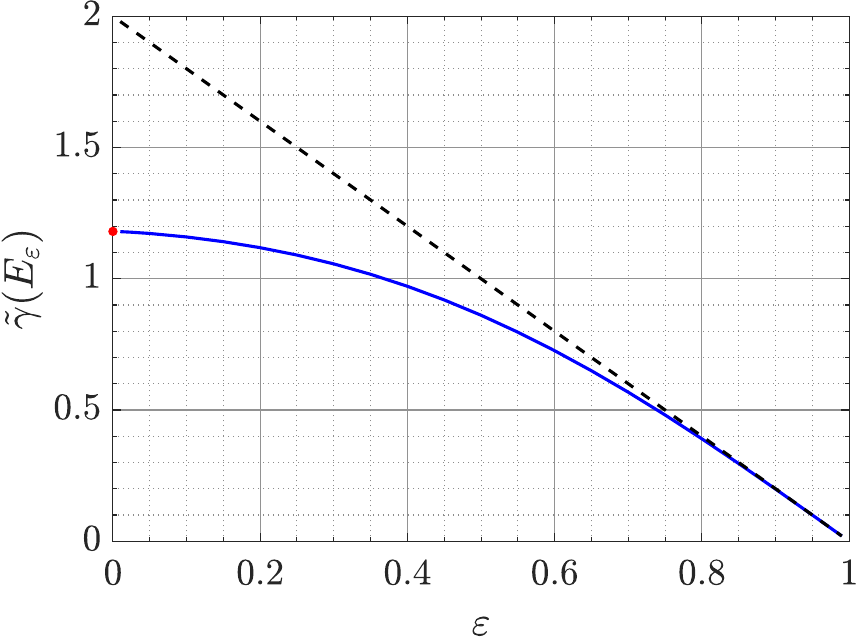}}
		\hfil
		\scalebox{0.457}{\includegraphics[trim=0 -1cm 0 0,clip]{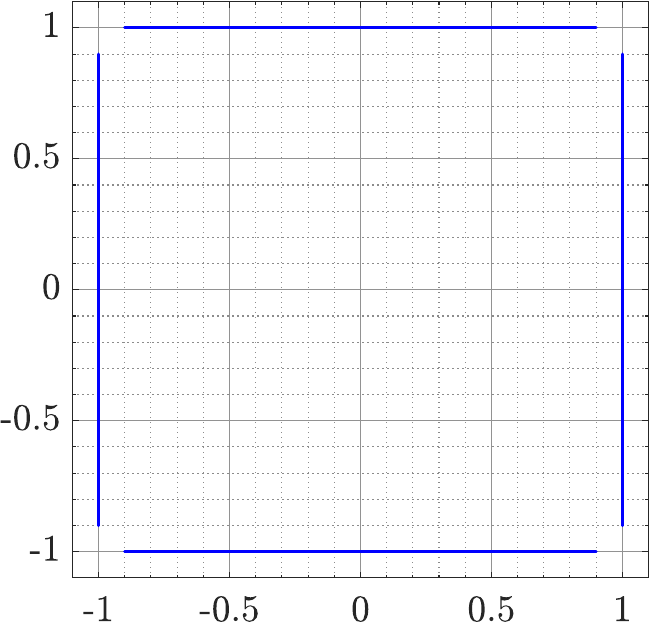}}
	}
	\caption{Graph of the approximate values $\tilde \gamma(E_{\varepsilon})$ of the analytic capacity for the union of rectilinear slits in Example~\ref{ex:sq-slits} as a function of $\varepsilon$ (left). A schematic of the compact set $E_\varepsilon$ for $\varepsilon=0.1$ (right). The red dot is the exact value of the analytic capacity of the square $[-1,1]\times[-1,1]$ and the dashed-line is the upper bound $2-2\varepsilon$ of $\gamma(E_\varepsilon)$.}
	\label{fig:sq-slits}
\end{figure}

\begin{example}\label{ex:sq-slits}
	We next consider the union of four equal rectilinear slits of length $2-2\varepsilon$, $0<\varepsilon<1$, such that these four slits make the square $[-1,1]\times[-1,1]$ when $\varepsilon=0$. A schematic of this configuration is shown in Figure~\ref{fig:sq-slits} (right) when $\varepsilon=0.1$. We denote these four slits by $F_{k,\varepsilon}$ with $k=1,2,3,4$. We define $E_\varepsilon=\cup_{k=1}^4F_{k,\varepsilon}$. 
\end{example}

We use our method to approximate $\gamma(E_\varepsilon)$ for $0.01\le\varepsilon\le0.99$ and the obtained numerical results are presented in Figure~\ref{fig:sq-slits} (left). As $\varepsilon\to0$, it is clear that the approximate values $\tilde\gamma(E_\varepsilon)$ approach 
\[
\gamma(E_0)= \frac{\Gamma^2(1/4)}{2\sqrt{\pi^3}} \approx 1.1803405990161,
\]
i.e., the value of the analytic capacity of the square $[-1,1]\times[-1,1]$. 
Furthermore, we always have 
\[
\tilde\gamma(E_\varepsilon)\le\sum_{k=1}^4\gamma(F_{k,\varepsilon})=2-2\varepsilon.
\]

\section{Concluding remarks}

This paper has shown how to use a numerical boundary integral equation method to quickly and accurately compute analytic capacity, an important conformal invariant. This quantity has been widely studied from a mostly theoretical perspective with several deep analytical results having been established~\cite{AVV,BelAhl,BelBook,Dav,Gar,HKV,Mur91,Mur94,PRY19,Tol03,Tol,Vit,YR13,YR18,Za}. Analytic capacity is intimately connected to the Ahlfors map and Szeg\"o kernel -- two fundamental objects in complex analysis --- and arises from the generalization of the Riemann map to multiply connected domains.

In our work, two particular classes of configurations were considered over which our calculations of analytic capacity were performed: compact sets bounded by smooth and piecewise smooth Jordan arcs, and domains consisting of a finite number of rectilinear slits. 
Throughout, we made connections with previous results; in particular, we have been able to corroborate the bounds found by Younsi \& Ransford \cite{YR13} for the analytic capacity for several compact sets they considered, and we were able to validate the conjectural subadditivity property of analytic capacity for numerous configurations. We were also able to validate numerically other exact results of analytic capacity, and illustrate several of its properties. Furthermore, the presented numerical results demonstrate that the analytic capacity $\gamma(E\cup F)$, roughly speaking, increases when the distance between $E$ and $F$ increases and decreases when this distance decreases.

Our work has been mainly numerical in approach, and the key to its success lies in the BIE scheme based on the Kerzman--Stein BIE~\cite{BelAhl,Ker-Tru,Murid1,Tru-Sze} and the FMM~\cite{Gre-Gim12,GR}. The method presented can be used for domains with smooth and piecewise smooth boundaries as well as for domains with many boundary components. We used the method also to approximate the analytic capacity for compact sets consisting of rectilinear slits. However, for the latter case, a preliminary conformal mapping step is required; this has also been shown to be expedient in other works~\cite{LSN17,NG18,Nvm}.

The MATLAB codes for the presented computations in this paper can be found at the link {\url{https://github.com/mmsnasser/ac}}.

\section*{Acknowledgements}
The authors would like to thank Nick Trefethen, Malik Younsi, and an anonymous reviewer for their valuable corrections, comments and suggestions, and for bringing several bibliographic items to our attention which greatly improved the presentation of this paper.

\begin{appendices}
\section{Computing a preimage domain for the rectilinear slit domain} \label{sec:App-pre}

Let $\Omega$ be a given multiply connected domain that is obtained by removing $m$ rectilinear slits $E_1,\ldots,E_m$ from the extended complex plane $\overline{\CC}$ such that the slit $E_j$ makes an angle $\theta_j$ with the positive $x$-axis, $j=1,\ldots,m$. 
In this appendix, we will summarize the iterative method from~\cite{NG18} for the construction of a preimage unbounded multiply connected domain $G$ bounded by $m$ smooth Jordan curves $\Gamma_1,\ldots,\Gamma_m$ as well as a conformal mapping $\zeta=\Phi(z)$ from $G$ onto $\Omega$.
This method has already been used in~\cite{LSN17,Nvm} for computing numerically the logarithmic and conformal capacities for rectilinear slit domains.

We assume that the boundary components $\Gamma_1,\ldots,\Gamma_m$ of the required unbounded domain $G$ are ellipses and are parametrized by
\begin{equation}\label{eq:ap-eta-j}
	\eta_j(t)=c_j+0.5a_j e^{\i \theta_j}(\cos t-r\sin t), \quad t\in J_j=[0,2\pi], 
\end{equation}
where $c_j$ is the center of the ellipse $\Gamma_j$ and $a_j$ is the length of its major axis, $j=1,\ldots,m$. The real parameter $r$ is the ratio between the length of the minor and major axes of these ellipses. We will choose its value such that $0<r\le 1$ where the domain $G$ is a circular for $r=1$. The value of $r$ is chosen to be $r=1$ when the slits are well separated and less than $1$ when the slits are close to each other (see~\cite{NG18} for details). 
Our objective here is to find the values of the parameters $c_j$ and $a_j$, $j=1,\ldots,m$, of the domain $G$ as well as a conformal mapping $\zeta=\Phi(z)$ from the domain $G$  in the $z$-plane onto the given unbounded rectilinear slit domain $\Omega$ in the $\zeta$-plane. 
With the normalization
\begin{equation}\label{eq:ap-cm-recs}
	\Phi(z)=z+O\left(\frac{1}{z}\right),
\end{equation} 
near infinity, such a conformal mapping is unique.

The conformal mapping $\Phi(z)$ can be computed using the following boundary integral equation method from~\cite{Nas-JMAA11}. We parametrize the whole boundary $\Gamma$ of the domain $G$ on the total parameter domain $J$ by
\begin{equation}\label{eq:ap-zeta}
	\eta(t)= \left\{ \begin{array}{l@{\hspace{0.5cm}}l}
		\eta_1(t),&t\in J_1, \\
		\quad\vdots & \\
		\eta_m(t),&t\in J_m.
	\end{array}
	\right.
\end{equation}
Then, we define a complex-valued function $A(t)$  on $J$ by
\begin{equation}\label{eq:ap-A}
	A(t)= e^{(\pi/2-\theta(t))\i},
\end{equation}
where $\theta(t)=\theta_j$ for $t\in J_j$, $j=1,\ldots,m$. With the functions $A(t)$ and $\eta(t)$, we define the kernels $N(s,t)$ and $M(s,t)$ for $(s,t)\in J\times J$ by
\begin{eqnarray}
	\label{eq:ap-N}
	N(s,t) &=&
	\frac{1}{\pi}\Im\left(\frac{A(t)}{A(s)}\frac{\eta'(t)}{\eta(t)-\eta(s)}\right),\\
	\label{eq:ap-M}
	M(s,t) &=&
	\frac{1}{\pi}\Re\left(\frac{A(t)}{A(s)}\frac{\eta'(t)}{\eta(t)-\eta(s)}\right).
\end{eqnarray}
The kernel $N(s,t)$, which is known as the generalized Neumann kernel, is continuous and the kernel $M(s,t)$ is singular~\cite{Weg-Nas}. The integral operators with the kernels $N(s,t)$ and $M(s,t)$ are then defined by
\[
\bN\mu(s) = \int_J N(s,t) \mu(t) dt, \quad
\bM\mu(s) = \int_J  M(s,t) \mu(t) dt, \quad s\in J.
\]
Let the function $\gamma$ be defined by
\begin{equation}\label{eq:ap-gam}
	\gamma(t)=\Im\left[e^{-\i\theta(t)}\eta(t)\right], \quad t\in J,
\end{equation}
let $\mu$ be the unique solution of the boundary integral equation with the generalized Neumann kernel
\begin{equation}\label{eq:ap-ie-g}
	(\bI-\bN)\mu=-\bM\gamma,
\end{equation}
and let the function $h$ be given by
\begin{equation}\label{eq:ap-h-g}
	h=\left(\bM\mu-(\bI-\bN)\gamma\right)/2.
\end{equation}
Then the function $\Psi$ with the boundary values
\begin{equation}\label{eq:ap-f-rec}
	\Psi(\eta(t))=\gamma(t)+h(t)+\i\mu(t)
\end{equation}
is analytic in $G$ with $\Psi(\infty)=0$. The values of $\Psi(z)$ for $z\in G$ can be computed by the Cauchy integral formula. Then, the values of the conformal mapping $\Phi(z)$ are given for $z\in G\cup\Gamma$ by
\begin{equation}\label{eq:ap-omega-app}
	\Phi(z)=z+\Psi(z).
\end{equation} 
For more details, see~\cite{Nas-JMAA11}.

The application of this method requires that the domain $G$ is known. However, in our case, the slit domain $\Omega$ is known and the domain $G$ is unknown and needs to be determined alongside the conformal mapping $\zeta=\Phi(z)$ from $G$ onto $\Omega$. This preimage domain $G$ as well as the conformal mapping $\Phi(z)$ will be computed using the following iterative method from~\cite{NG18} which generates a sequence of multiply connected domains $G^{(0)},G^{(1)},G^{(2)},\ldots,$ that converges numerically to the required preimage domain $G$. 
Let $L_j=|E_j|$ be the length of the slit $E_{j}$ and let $\beta_j$ be its center, $j=1,\ldots,m$. 
In the iteration step $i=0,1,2,\ldots$, we assume that $G^{(i)}$ is an unbounded multiply connected domain bounded by the $m$ ellipses $\Gamma^{(i)}_1,\ldots,\Gamma^{(i)}_m$ parametrized by
\begin{equation}\label{eq:ap-eta-i}
	\zeta^{(i)}_j(t)=c^{(i)}_j+0.5a^{(i)}_je^{\i \theta_j}(\cos t-r\sin t), \quad 0\le t\le 2\pi, \quad j=1,\ldots,m.
\end{equation}
The parameters $c^{(i)}_j$ and $a^{(i)}_j$ are computed as follows:
\begin{enumerate}
	\item Set
	\[
	c^{(0)}_j=\beta_j, \quad a^{(0)}_j=(1-0.5r)L_j, \quad j=1,\ldots,m.
	\]
	\item For $i=1,2,3,\ldots,$
	\begin{itemize}
		\item Compute the conformal mapping from the preimage domain $G^{(i-1)}$ in the $z$-plane onto the canonical rectilinear slit domain $\Omega^{(i)}$ obtained by removing $m$ rectilinear slits $E^{(i)}_{j}$, $j=1,\ldots,m$, from the $\zeta$-plane (using the method presented in~\eqref{eq:ap-gam}--\eqref{eq:ap-omega-app} above). Let $L^{(i)}_j$ denote the length of the slit $E^{(i)}_{j}$ and let $\beta^{(i)}_j$ denote its center.
		\item Define	
		\[
		c^{(i)}_j = c^{(i-1)}_j-(\beta^{(i)}_j-\beta_j), \quad 
		a^{(i)}_j = a^{(i-1)}_j-(1-0.5r)(L^{(i)}_j -L_j), \quad j=1,\ldots,m.
		\]				
	\end{itemize}
	\item Stop the iteration if 
	\[
	\frac{1}{2m}\sum_{j=1}^{m}\left(|\beta^{(i)}_j - \beta_j|+|L^{(i)}_j -L|\right)<\varepsilon \quad{\rm or}\quad i>{\tt Max}
	\]
	where $\varepsilon$ is a given tolerance and ${\tt Max}$ is the maximum number of iterations allowed.		
\end{enumerate}

The above iterative method generates sequences of parameters $c^{(i)}_j$ and $a^{(i)}_j$ that converge numerically to $c_j$ and $a_j$, respectively, and then the boundary components of the preimage domain $G$ are parametrized by~\eqref{eq:ap-eta-j}. In our numerical implementations, we used $\varepsilon=10^{-13}$ and ${\tt Max}=100$.

\begin{figure}[htb] %
	\centerline{
		\scalebox{0.4}{\includegraphics[trim=0 0 0 0,clip]{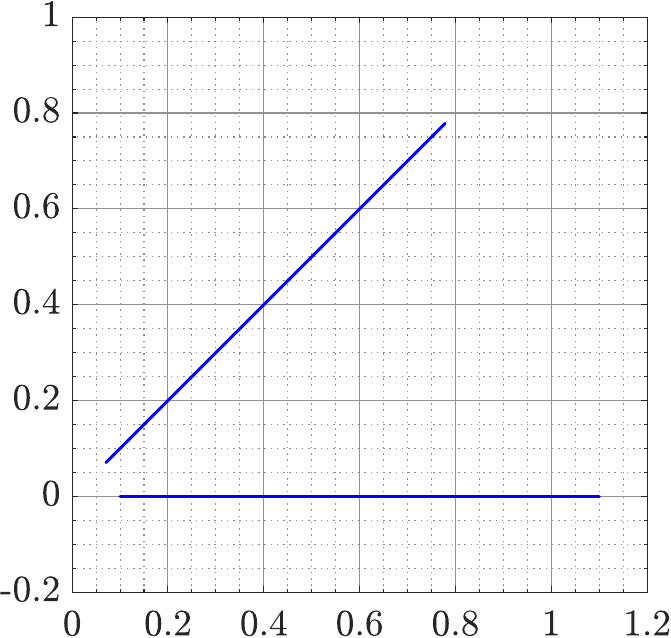}}
		\hfil
		\scalebox{0.4}{\includegraphics[trim=0 0 0 0,clip]{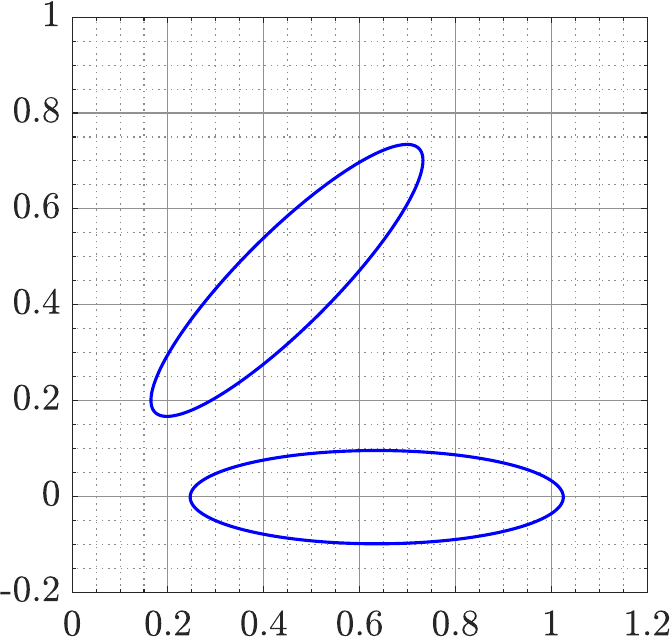}}
	}
	\caption{The given slit domain $\Omega$ in Example~\ref{ex:flap} (for $\varepsilon=0.25$) and the computed preimage domain $G$ (with $r=0.25$).}
	\label{fig:ap-2s}
\end{figure}

\begin{figure}[htb] %
	\centerline{
		\scalebox{0.4}{\includegraphics[trim=0 0 0 0,clip]{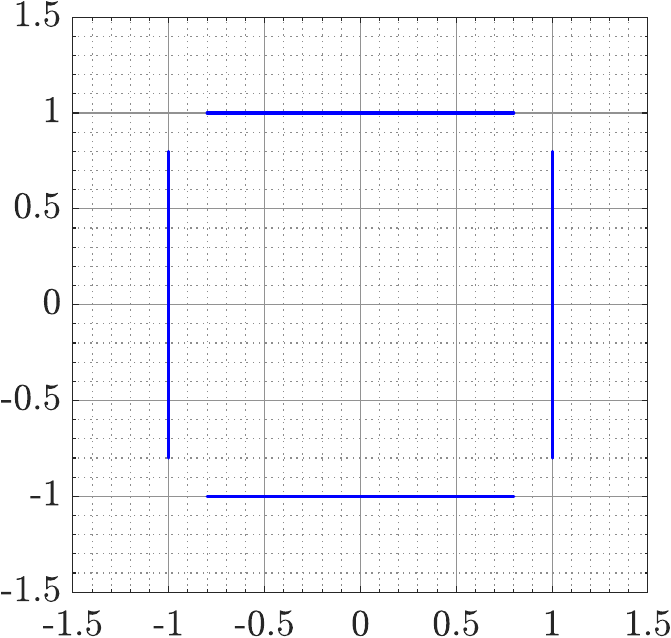}}
		\hfil
		\scalebox{0.4}{\includegraphics[trim=0 0 0 0,clip]{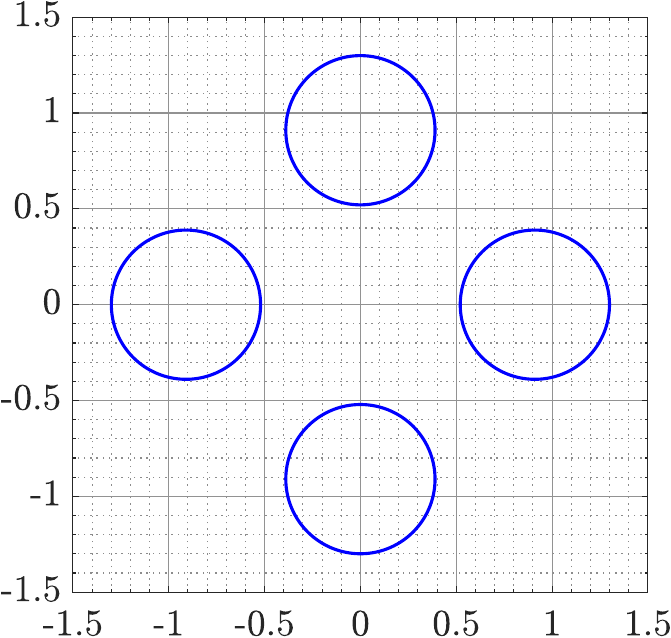}}
	}
	\caption{The given slit domain $\Omega$ in Example~\ref{ex:sq-slits} (for $\varepsilon=0.2$) and the computed preimage domain $G$ (with $r=1$).}
	\label{fig:ap-4s}
\end{figure}

It is clear that each iteration of the above method requires solving the integral equation with the generalized Neumann kernel~(\ref{eq:ap-ie-g}) and  computing the function $h$ in~(\ref{eq:ap-h-g}) which can be done with the {\sc Matlab} function \verb|fbie| presented in~\cite{Nas-ETNA}. 
In our numerical computations, the value of $n$ as well as the values of the other parameters in \verb|fbie| are chosen to be the same as those used in the method described in Section~\ref{sec:nm} for computing the analytic capacity. The given slit domain $\Omega$ and the computed preimage domain $G$ for Examples~\ref{ex:flap} and~\ref{ex:sq-slits} are presented in Figures~\ref{fig:ap-2s} and~\ref{fig:ap-4s}, respectively.

\end{appendices}


\end{document}